\documentclass[10pt,a4paper]{article}
\usepackage[latin1]{inputenc}
\usepackage{amsmath}
\usepackage{amsfonts}
\usepackage{amssymb}
\usepackage{amsthm}

\newtheorem{thm}{Theorem}[section]
\newtheorem{prop}{Proposition}[section]
\newtheorem{lem}{Lemma}[section]
\newtheorem{eg}{Example}[section]
\newtheorem*{condition}{Condition (C)}

\theoremstyle{remark}
\newtheorem{rmk}{Remark}[section]

\numberwithin{equation}{section}

\def \div{\hbox{div}}
\def \tr{\hbox{tr}}
\def \Hess{\hbox{Hess}}
\def \diag{\hbox{diag}}

\newcommand{\R}{\mathbb{R}}
\newcommand{\N}{\mathbb{N}}

\title{\textbf{Schr\"{o}dinger Operator:\\ Heat Kernel and Its Applications}}
\author{Sheng-Ya Feng}
\date{}

\begin{document}
\maketitle

\begin{abstract}
\noindent
In this paper, we study the geometry associated with Schr\"{o}dinger operator via Hamiltonian and Lagrangian formalism. Making use of a multiplier technique, we construct the heat kernel with the coefficient matrices of the operator both diagonal and non-diagonal. For applications, we compute the heat kernel of a Schr\"{o}dinger operator with terms of lower order, and obtain a globally closed solution to a matrix Riccati equations as a by-product. Besides, we finally recover and generalise several classical results on some celebrated operators.

\bigskip
\hspace{-15pt}\textbf{Key Words:} Hamiltonian system, Hamilton-Jacobi equation, transport equation, matrix Riccati equation, Schr\"{o}dinger operator
\end{abstract}

\bigskip
\section{Introduction}
We first introduce a second order differential operator with quadratic potentials
$$
T=-\div(A\nabla) + \langle Bx,x \rangle +\langle Cx,\nabla\rangle.
$$
with $A$, $B$ and $C$ matrices. From now on, we call the fundamental solution of the operator $\partial_{t}+T$ the \textit{heat kernel} of $T$. Let us recall some well-known facts for $B=0$, when $T$ becomes
$$
H=-\div(A\nabla) +\langle Cx,\nabla\rangle.
$$
Kolmogorov \cite{K} considers the following equation 
$$
(\partial_{t}-\partial_{x_{1}}^{2}-x_{1}\partial_{x_{2}})u=0, \hspace{15pt}(x_{1}, x_{2}, t)\in \R\times\R\times\R^{+}
$$
to describe the probability density of a system with $2n$ degree of freedom and obtains an explicit fundamental solution by Fourier transform. H\"{o}rmander \cite{H} uses the same method to construct heat kernel for $H$ under a condition imposed on $A$ and $B$ which is equivalent to the hypoellipticity of $\partial_{t} +H$. Beals \cite{B1} sketches a method to find heat kernel for $H$ with $A\geq 0$ via a probabilistic ansatz.

The study of the generalised Hermite operator $ L= -\varDelta + \langle Bx, x \rangle$ is of independent interest. It takes Hermite operator and anti-Hermite operator as its typical cases. Hermite operator $ L_{H} = -\varDelta + \lvert x \rvert ^{2}$ arises from harmonic oscillator and has been studied for quite some time (cf. \cite{B3}, \cite{GJ}), while anti-Hermite operator $ L_{GH} = -\varDelta - \lvert x \rvert ^{2}$ arises from anti-harmonic oscillator discussed in \cite{RS}. To the best knowledge of this author, the geometry induced by anti-harmonic oscillator was seldom studied. \cite{CF} makes some effort in this direction. They study the geometry of generalised Hermite operator $L_{GH}= -\varDelta + \langle Bx, x\rangle$ with $B$ any real matrix by characterizing the behaviour of geodesics when spatial dimension $n$ equals 2. 

Our interests concentrate on the case $C=0$, and we study Schr\"{o}dinger operator
\begin{equation}
L_{S}=-\div(A\nabla)+\langle Bx,x \rangle
\end{equation}
where $ A $ is  a symmetric positive definite $ n\times n$ real matrix, $ B $ is a $ n\times n$ real matrix commutative with $ A $, i.e. $ AB=BA $, $ \div $, $\nabla$, and $ \langle\cdot , \cdot\rangle $ denote respectively divergence, gradient and Euclidean inner product. In section 2, we quantitatively study the associated Hamiltonian system for any dimension, from which we conclude that the singularities are hyperplanes in phase space. This work uniformly generalises the results for $ B $ positive definite or low dimensional space (cf.\cite{CF}, \cite{F} or \cite{FD}). Moreover, we formally characterize three important objects, namely geodesic, energy and action function of the Hamiltonian system, where \textit{geodesic} formally means $x$-component of the solution for Hamiltonian system. All these quantities are given in closed forms, which will play a crucial role in constructing explicit heat kernel of Schr\"{o}dinger operator in section 3.

A common way to compute the heat kernel of Hermite operator is to use eigenfunction expansion (cf. \cite{T}). In recent, \cite{CCT} studies Hamiltonian system qualitatively from the view point of conservation law of energy, and obtains the heat kernel formulae with $B$ a diagonally positive definite matrix. \cite{F} and \cite{FD} generalise the results in terms of spectral calculus for $B$ any positive semi-definite matrix. In section 3, we first use the obtained action function and a multiplier technique to construct the heat kernel of Schr\"{o}dinger operator $L_{S}$ when the coefficient matrices $A$ and $B$ are diagonal. It is worth mentioning that the heat kernel has slightly different properties as the normal one does, primly because that the Schr\"{o}dinger operator under consideration is not linear in $x-$variables. For this reason, we address the fundamental solution of $\partial_{t}+L_{S}$ or $\partial_{t}+L$ as \textit{generalised heat kernel}. We close section 3 with the computation of heat kernel for $A$ and $B$ non-diagonal case.

The heat kernel has significance in two areas of applications. In section 4, we first apply it to obtain the heat kernel for Schr\"{o}dinger operator with terms of lower order
\begin{equation}
L=-\div(A\nabla) + \langle Bx,x \rangle + \langle f,\nabla \rangle + \langle g,x \rangle+h
\end{equation}
where real matrix $ A $ is positive definite and commutative with $ B $, $f$ and $g$ are vectors, and $h$ is a real number. The heat kernel of $ L $ has an ansatz
\begin{equation}
K(x,x^{0};t)= W(t)\exp\{\langle\alpha(t)x,x\rangle+\langle\beta(t)x,x^{0}\rangle+\langle\gamma(t)x^{0},x^{0}\rangle + \langle\mu,x\rangle+\langle\nu,x^{0}\rangle\},\label{1.3}
\end{equation}
where $\alpha$, $\beta$, $\gamma$ are expected to be symmetric $ n\times n$ real matrices, $\mu$, $\nu$ to be vectors, and we deduce a system of matrix and scalar differential equations as in \cite{B1}
\begin{align}
&\dot{\alpha} = 4\alpha A \alpha - \frac{B+B^{t}}{2}\label{1.4}\\
&\dot{\beta} = 4\beta A \alpha \label{1.5}\\
&\dot{\gamma} = \beta A \beta \label{1.6}\\
&\dot{\mu} = 4\alpha A \mu - 2 \alpha f - g \label{1.7}\\
&\dot{\nu} =2 \beta A \mu -\beta f \label{1.8}\\
&W^{-1}\dot{W} = 2 \tr(A \alpha) + \langle A \mu, \mu\rangle - \langle f,\mu\rangle -h \label{1.9}
\end{align}
where the dot denotes $\frac{d}{dt}$, and $B^{t}$ denotes transpose of $B$. The main difficulty is to solve the matrix Riccati equation (\ref{1.4}), which is an equation of fundamental importance in control theory \cite{AFIJ}. Fortunately, the heat kernel of $L_{S}$ provide us a globally closed solution of matrix Riccati equation (\ref{1.4}), and a condition to identify the solution of the scalar differential equation (\ref{1.9}). Then other equations and hence the heat kernel of $L$ can be explicitly computed. Last section is devoted to the second areas of applications. We will recover and generalise several classical results on some celebrated operators, including Laplacian, Hermite operator and Ornstein-Uhlenbeck operator on weighted space.

\bigskip
\section{Hamiltonian system associated with $L_{S}$ }
In this section, we consider Hamiltonian system associated with Schr\"{o}dinger operator
$$
L_{S} = -\div \left( A \nabla\right)  + \langle Bx, x \rangle
$$
with $ A $ and $ B $ commutative.

Geodesics, energy and Hamilton-Jacobi action function are three significant objects in Hamilton-Jacobi theory and are of their own interest. We study them one by one in the following  subsections.

\medskip
\subsection{Geodesics}
The Hamiltonian function of $L_{S}$ is defined as its full symbol 
$$
H_{S}= -\langle A\xi, \xi\rangle + \langle Bx, x\rangle
$$ 
and the associated Hamiltonian system is 
\begin{equation}
\left\{\begin{aligned}
& \dot{x}=\frac{\partial H_{S}}{\partial \xi}= -A\xi- A^{t}\xi= -2A\xi\\
& \dot{\xi}= -\frac{\partial H_{S}}{\partial x}= -Bx-B^{t}x= -\left( B+B^{t}\right) x
\end{aligned}\right..\label{3.1}
\end{equation}
Denoting $D:= 2A\left( B+B^{t}\right) $, one has 
$$
\ddot{x}= -2A\dot{\xi}= 2A(B+B^{t})x= Dx,
$$
and $D$ is symmetric following from that $A$ is commutative with $B$.
The geodesics $ x(s) $ between $ x^{0} $ and $ x $ in $ \mathbb{R}^{n} $ satisfy the boundary value problem
\begin{equation}
\left\{\begin{aligned}
& \ddot{x}=Dx\\
& x(0)=x_{0},\hspace{5pt} x(t)=x 
\end{aligned}\right..\label{3.2}
\end{equation}

We start with the case when both $A$ and $B$ are diagonal matrices, and write them as follows
\begin{equation}
A=\varLambda^{a}=
\begin{bmatrix}
a_{1}^{2} & & \\
 & \ddots & \\
 & & a_{n}^{2}
\end{bmatrix}\label{3.3}
\end{equation}

\begin{equation}
B=\varLambda^{b}=
\begin{bmatrix}
b_{1}^{2} & & & & & \\
 & \ddots & & & & \\
 & & b_{m}^{2} & & & \\
 & & & -b_{m+1}^{2} & & \\
 & & & & \ddots & \\
 & & & & & -b_{n}^{2}
\end{bmatrix}\label{3.4}
\end{equation}

\begin{equation}
AB=\varLambda^{a}\varLambda^{b}=
\begin{bmatrix}
a_{1}^{2}b_{1}^{2} & & & & & \\
 & \ddots & & & & \\
 & & a_{m}^{2}b_{m}^{2} & & & \\
 & & & -a_{m+1}^{2}b_{m+1}^{2} & & \\
 & & & & \ddots & \\
 & & & & & -a_{n}^{2}b_{n}^{2}
\end{bmatrix}
\end{equation}
where $a_{j}>0$, $b_{j}>0$ for $j\in \{1,\cdots, n\}$ satisfy the condition that:

\bigskip
\begin{condition}
For $i \neq k$ and $1\leqslant i, k \leqslant m$ or $m+1\leqslant i, k \leqslant n$, $a_{i}^{2}b_{i}^{2}\neq a_{k}^{2}b_{k}^{2}$.
\end{condition}

\bigskip
Putting $\varLambda_{1}^{a}=\diag \{a_{j}^{2}\}_{j=1}^{m}$, $\varLambda_{2}^{a}=\diag \{a_{j}^{2}\}_{j=m+1}^{n}$, $\varLambda_{1}^{b}=\diag \{b_{j}^{2}\}_{j=1}^{m}$, and $\varLambda_{2}^{b}=\diag \{-b_{j}^{2}\}_{j=m+1}^{n}$, one has 
\begin{equation}
A=
\begin{bmatrix}
\varLambda_{1}^{a} & \\
 & \varLambda_{2}^{a}
\end{bmatrix},\label{3.6}
\end{equation}

\begin{equation}
B=
\begin{bmatrix}
\varLambda_{1}^{b} & \\
 & \varLambda_{2}^{b}
\end{bmatrix},\label{3.7}
\end{equation}

and
 
\begin{equation}
D=4AB=4
\begin{bmatrix}
\varLambda_{1}^{a}\varLambda_{1}^{b} & \\
 & \varLambda_{2}^{a}\varLambda_{2}^{b}
\end{bmatrix}.
\end{equation}

The solution of the linear system $\ddot{x}=Dx$ is a combination of radical solutions $\{e^{2a_{j}b_{j}s}\}_{j=1}^{m}$, $\{e^{-2a_{j}b_{j}s}\}_{j=1}^{m}$, $\{\cos(2a_{j}b_{j}s)\}_{j=m+1}^{n}$, and $\{\sin(2a_{j}b_{j}s)\}_{j=m+1}^{n}$. We write the coefficients in both block and component forms

\begin{align*}
C_{1}&=
\begin{bmatrix}
C_{11}\\
C_{21}
\end{bmatrix},
&C_{2}&=
\begin{bmatrix}
C_{12}\\
C_{22}
\end{bmatrix},\\
C_{3}&=
\begin{bmatrix}
C_{13}\\
C_{23}
\end{bmatrix},
&C_{4}&=
\begin{bmatrix}
C_{14}\\
C_{24}
\end{bmatrix},
\end{align*}

where

\begin{align*}
C_{11}&=
\begin{bmatrix}
c_{11} & \cdots & c_{1m}\\
\vdots & \ddots & \vdots\\
c_{m1} & \cdots & c_{mm}
\end{bmatrix},
&C_{21}&=
\begin{bmatrix}
c_{m+1,1} & \cdots & c_{m+1,m}\\
\vdots & \ddots & \vdots\\
c_{n1} & \cdots & c_{nm}
\end{bmatrix},\\
C_{12}&=
\begin{bmatrix}
c_{1,m+1} & \cdots & c_{1n}\\
\vdots & \ddots & \vdots\\
c_{m,m+1} & \cdots & c_{mn}
\end{bmatrix},
&C_{22}&=
\begin{bmatrix}
c_{m+1,m+1} & \cdots & c_{m+1,n}\\
\vdots & \ddots & \vdots\\
c_{n1} & \cdots & c_{nn}
\end{bmatrix},\\
C_{13}&=
\begin{bmatrix}
c_{1,n+1} & \cdots & c_{1,n+m}\\
\vdots & \ddots & \vdots\\
c_{m,n+1} & \cdots & c_{m,n+m}
\end{bmatrix},
&C_{23}&=
\begin{bmatrix}
c_{m+1,n+1} & \cdots & c_{m+1,n+m}\\
\vdots & \ddots & \vdots\\
c_{n,n+1} & \cdots & c_{n,n+m}
\end{bmatrix},\\
C_{14}&=
\begin{bmatrix}
c_{1,n+m+1} & \cdots & c_{1,2n}\\
\vdots & \ddots & \vdots\\
c_{m,n+m+1} & \cdots & c_{m,2n}
\end{bmatrix},
&C_{24}&=
\begin{bmatrix}
c_{m+1,n+m+1} & \cdots & c_{m+1,2n}\\
\vdots & \ddots & \vdots\\
c_{n,n+m+1} & \cdots & c_{n,2n}
\end{bmatrix}.
\end{align*}

Accordingly, we first write the solution vectors as 
\begin{align*}
& x_{1}(s)= (e^{2a_{1}b_{1}s},\cdots, e^{2a_{m}b_{m}s})^{t},\\
& x_{2}(s)=(\cos(2a_{m+1}b_{m+1}s),\cdots,\cos(2a_{n}b_{n}s))^{t},\\
& x_{3}(s)= (e^{-2a_{1}b_{1}s},\cdots, e^{-2a_{m}b_{m}s})^{t},\\
& x_{4}(s)=(\sin(2a_{m+1}b_{m+1}s),\cdots,\sin(2a_{n}b_{n}s))^{t},
\end{align*} 

then 

$$
x(s)=
\begin{bmatrix}
C_{11} & C_{12} & C_{13} & C_{14}\\
C_{21} & C_{22} & C_{23} & C_{24}
\end{bmatrix}
\begin{bmatrix}
x_{1}(s)\\
x_{2}(s)\\
x_{3}(s)\\
x_{4}(s)
\end{bmatrix},
$$

\begin{align*}
\ddot{x}(s)&=4
\begin{bmatrix}
C_{11} & C_{12} & C_{13} & C_{14}\\
C_{21} & C_{22} & C_{23} & C_{24}
\end{bmatrix}
\begin{bmatrix}
\varLambda_{1}^{a}\varLambda_{1}^{b} & & & \\
& \varLambda_{2}^{a}\varLambda_{2}^{b} & & \\
& & \varLambda_{1}^{a}\varLambda_{1}^{b}  & \\
& & & \varLambda_{2}^{a}\varLambda_{2}^{b}\\
\end{bmatrix}
\begin{bmatrix}
x_{1}(s)\\
x_{2}(s)\\
x_{3}(s)\\
x_{4}(s)
\end{bmatrix}\\
&=4
\begin{bmatrix}
C_{11}\varLambda_{1}^{a}\varLambda_{1}^{b} & C_{12}\varLambda_{2}^{a}\varLambda_{2}^{b} & C_{13}\varLambda_{1}^{a}\varLambda_{1}^{b} & C_{14}\varLambda_{2}^{a}\varLambda_{2}^{b}\\
C_{21}\varLambda_{1}^{a}\varLambda_{1}^{b} & C_{22}\varLambda_{2}^{a}\varLambda_{2}^{b} & C_{23}\varLambda_{1}^{a}\varLambda_{1}^{b} & C_{24}\varLambda_{2}^{a}\varLambda_{2}^{b}
\end{bmatrix}
\begin{bmatrix}
x_{1}(s)\\
x_{2}(s)\\
x_{3}(s)\\
x_{4}(s)
\end{bmatrix},
\end{align*}

\begin{align*}
Dx(s) &=4
\begin{bmatrix}
\varLambda_{1}^{a}\varLambda_{1}^{b} & \\
& \varLambda_{2}^{a}\varLambda_{2}^{b}
\end{bmatrix}
\begin{bmatrix}
C_{11} & C_{12} & C_{13} & C_{14}\\
C_{21} & C_{22} & C_{23} & C_{24}
\end{bmatrix}
\begin{bmatrix}
x_{1}(s)\\
x_{2}(s)\\
x_{3}(s)\\
x_{4}(s)
\end{bmatrix}\\
&=4
\begin{bmatrix}
\varLambda_{1}^{a}\varLambda_{1}^{b}C_{11} & \varLambda_{1}^{a}\varLambda_{1}^{b}C_{12} & \varLambda_{1}^{a}\varLambda_{1}^{b}C_{13} & \varLambda_{1}^{a}\varLambda_{1}^{b}C_{14}\\
\varLambda_{2}^{a}\varLambda_{2}^{b}C_{21} & \varLambda_{2}^{a}\varLambda_{2}^{b}C_{22} & \varLambda_{2}^{a}\varLambda_{2}^{b}C_{23} & \varLambda_{2}^{a}\varLambda_{2}^{b}C_{24}
\end{bmatrix}
\begin{bmatrix}
x_{1}(s)\\
x_{2}(s)\\
x_{3}(s)\\
x_{4}(s)
\end{bmatrix}.
\end{align*}

Noting $\ddot{x}(s)=Dx(s)$, and the condition (C) implies
\begin{equation}
\begin{aligned}
C_{11}&=\diag\{c_{jj}\}_{j=1}^{m}, &C_{13}&=\diag\{c_{j,n+j}\}_{j=1}^{m},\\
C_{22} &=\diag\{c_{jj}\}_{j=m+1}^{n}, &C_{24}&=\diag\{c_{j,n+j}\}_{j=m+1}^{n}.
\end{aligned}\label{3.9}
\end{equation}
Similarly, that $a_{j}^{2}b_{j}^{2}>0$ for $j\in \overline{1,n}$ implies
$$
C_{12}=C_{21}=C_{14}=C_{23}=0.
$$
Thus, 
$$x(s)=
\begin{bmatrix}
C_{11} & 0 & C_{13} & 0\\
0 & C_{22} & 0 & C_{24}
\end{bmatrix}
\begin{bmatrix}
x_{1}(s)\\
x_{2}(s)\\
x_{3}(s)\\
x_{4}(s)
\end{bmatrix}
$$
where $C_{11}$, $C_{22}$, $C_{13}$, $C_{24}$ are diagonal matrices commutative with $\varLambda_{j}^{a}$ and $\varLambda_{j}^{b}$ for $j=1, 2$.

Next, the boundary condition in (\ref{3.2}) will establish $C_{ij}$'s. As before, we introduce some notations.
\begin{align*}
x^{0}&=(x_{1}^{(0)}, \cdots, x_{n}^{(0)})^{t}, &x &=(x_{1}^{(1)}, \cdots, x_{n}^{(1)})^{t},\\
x_{1}^{0} &=(x_{1}^{(0)}, \cdots, x_{m}^{(0)})^{t}, &x_{2}^{0}&=(x_{m+1}^{(0)}, \cdots, x_{n}^{(0)})^{t},\\
x_{1} &=(x_{1}^{(1)}, \cdots, x_{m}^{(1)})^{t}, &x_{2} &=(x_{m+1}^{(1)}, \cdots, x_{n}^{(1)})^{t},\\
\widetilde{C}_{11}&=(c_{11},\cdots,c_{mm})^{t}, &\widetilde{C}_{13}&=(c_{1,n+1},\cdots,c_{m,n+m})^{t},\\
\widetilde{C}_{22}&=(c_{m+1,m+1},\cdots,c_{nn})^{t}, &\widetilde{C}_{24}&=(c_{m+1,m+n+1},\cdots,c_{n,2n})^{t}.
\end{align*}
Given a non-singular matrix $M$, we define $\frac{N}{M}:=M^{-1}N$.

By boundary condition in (\ref{3.2}), $\widetilde{C}_{11}$, $\widetilde{C}_{22}$, $\widetilde{C}_{13}$ and $\widetilde{C}_{24}$ satisfy the following linear equations
$$
\begin{bmatrix}
I_{m} & 0 & I_{m} & 0\\
0 & I_{n-m} & 0 & 0\\
e^{2t\sqrt{\varLambda_{1}^{a}\varLambda_{1}^{b}}} & 0 & e^{-2t\sqrt{\varLambda_{1}^{a}\varLambda_{1}^{b}}} & 0 \\
0 & \cos( 2t\sqrt{-\varLambda_{2}^{a}\varLambda_{2}^{b}})  & 0 & \sin( 2t\sqrt{-\varLambda_{2}^{a}\varLambda_{2}^{b}}) 
\end{bmatrix}
\begin{bmatrix}
\widetilde{C}_{11}\\
\widetilde{C}_{22}\\
\widetilde{C}_{13}\\
\widetilde{C}_{24}
\end{bmatrix}
=
\begin{bmatrix}
x_{1}^{0}\\
x_{2}^{0}\\
x_{1}\\
x_{2}
\end{bmatrix}.
$$ 
To move on, we make assumption (*) that $\sin\left( 2t\sqrt{-\varLambda_{2}^{a}\varLambda_{2}^{b}}\right) $ is non-singular. Indeed, the region for $\sin\left( 2t\sqrt{-\varLambda_{2}^{a}\varLambda_{2}^{b}}\right) $ singular consists of countably many hyperplanes in $(x,t)-$space, thus, it has no contribution to the Hamilton-Jacobi action function. In section 3, we will give a formal remark about this point.
\begin{align*}
&\quad
\begin{bmatrix}
I_{m} & 0 & I_{m} & 0 & \vdots & x_{1}^{0}\\
0 & I_{n-m} & 0 & 0 & \vdots & x_{2}^{0}\\
e^{2t\sqrt{\varLambda_{1}^{a}\varLambda_{1}^{b}}} & 0 & e^{-2t\sqrt{\varLambda_{1}^{a}\varLambda_{1}^{b}}} & 0 & \vdots & x_{1}\\
0 & \cos(2t\sqrt{-\varLambda_{2}^{a}\varLambda_{2}^{b}}) & 0 & \sin(2t\sqrt{-\varLambda_{2}^{a}\varLambda_{2}^{b}}) & \vdots & x_{2}
\end{bmatrix}\\
&\rightarrow
\begin{bmatrix}
I_{m} & 0 & I_{m} & 0 & \vdots & x_{1}^{0}\\
0 & I_{n-m} & 0 & 0 & \vdots & x_{2}^{0}\\
0 & 0 & \frac{1-e^{4t\sqrt{\varLambda_{1}^{a}\varLambda_{1}^{b}}}}{e^{2t\sqrt{\varLambda_{1}^{a}\varLambda_{1}^{b}}}} & 0 & \vdots & x_{1}-e^{2t\sqrt{\varLambda_{1}^{a}\varLambda_{1}^{b}}}x_{1}^{0}\\
0 & 0 & 0 & \sin(2t\sqrt{-\varLambda_{2}^{a}\varLambda_{2}^{b}}) & \vdots & x_{2}-\cos(2t\sqrt{-\varLambda_{2}^{a}\varLambda_{2}^{b}})x_{2}^{0}
\end{bmatrix}\\
&\rightarrow
\begin{bmatrix}
I_{m} & 0 & I_{m} & 0 & \vdots & x_{1}^{0}\\
0 & I_{n-m} & 0 & 0 & \vdots & x_{2}^{0}\\
0 & 0 & -I_{m} & 0 & \vdots & \frac{e^{2t\sqrt{\varLambda_{1}^{a}\varLambda_{1}^{b}}}x_{1}-e^{4t\sqrt{\varLambda_{1}^{a}\varLambda_{1}^{b}}}x_{1}^{0}}{e^{4t\sqrt{\varLambda_{1}^{a}\varLambda_{1}^{b}}}-1}\\
0 & 0 & 0 & \sin(2t\sqrt{-\varLambda_{2}^{a}\varLambda_{2}^{b}}) & \vdots & x_{2}-\cos(2t\sqrt{-\varLambda_{2}^{a}\varLambda_{2}^{b}})x_{2}^{0}
\end{bmatrix}\\
&\xrightarrow{(*)}
\begin{bmatrix}
I_{m} & 0 & 0 & 0 & \vdots & \frac{e^{2t\sqrt{\varLambda_{1}^{a}\varLambda_{1}^{b}}}}{e^{4t\sqrt{\varLambda_{1}^{a}\varLambda_{1}^{b}}}-1}x_{1}-\frac{1}{e^{4t\sqrt{\varLambda_{1}^{a}\varLambda_{1}^{b}}}-1}x_{1}^{0}\\
0 & I_{n-m} & 0 & 0 & \vdots & x_{2}^{0}\\
0 & 0 & I_{m} & 0 & \vdots & -\frac{e^{2t\sqrt{\varLambda_{1}^{a}\varLambda_{1}^{b}}}}{e^{4t\sqrt{\varLambda_{1}^{a}\varLambda_{1}^{b}}}-1}x_{1}+\frac{e^{4t\sqrt{\varLambda_{1}^{a}\varLambda_{1}^{b}}}}{e^{4t\sqrt{\varLambda_{1}^{a}\varLambda_{1}^{b}}}-1}x_{1}^{0}\\
0 & 0 & 0 & I_{n-m} & \vdots & \csc(2t\sqrt{-\varLambda_{2}^{a}\varLambda_{2}^{b}})x_{2}-\cot(2t\sqrt{-\varLambda_{2}^{a}\varLambda_{2}^{b}})x_{2}^{0}
\end{bmatrix}.
\end{align*}
We read off the solution 
\begin{equation}
\begin{split}
\widetilde{C}_{11}&=\frac{e^{2t\sqrt{\varLambda_{1}^{a}\varLambda_{1}^{b}}}}{e^{4t\sqrt{\varLambda_{1}^{a}\varLambda_{1}^{b}}}-1}x_{1}-\frac{1}{e^{4t\sqrt{\varLambda_{1}^{a}\varLambda_{1}^{b}}}-1}x_{1}^{0},\\ \widetilde{C}_{22} &=x_{2}^{0},\\
\widetilde{C}_{13}&=-\frac{e^{2t\sqrt{\varLambda_{1}^{a}\varLambda_{1}^{b}}}}{e^{4t\sqrt{\varLambda_{1}^{a}\varLambda_{1}^{b}}}-1}x_{1}+\frac{e^{4t\sqrt{\varLambda_{1}^{a}\varLambda_{1}^{b}}}}{e^{4t\sqrt{\varLambda_{1}^{a}\varLambda_{1}^{b}}}-1}x_{1}^{0},\\
\widetilde{C}_{24}&=\csc\left( 2t\sqrt{-\varLambda_{2}^{a}\varLambda_{2}^{b}}\right) x_{2}-\cot\left( 2t\sqrt{-\varLambda_{2}^{a}\varLambda_{2}^{b}}\right) x_{2}^{0},
\end{split}\label{3.10}
\end{equation}
and we may recover $C_{ij}$'s in (\ref{3.9}) from $\widetilde{C}_{ij}$'s in (\ref{3.10}) with
\begin{equation}
\begin{aligned}
c_{jj}&=\left\langle \widetilde{C}_{11},e_{j}^{m}\right\rangle ,  &c_{j,n+j}&=\left\langle \widetilde{C}_{13},e_{j}^{m}\right\rangle ,  &j &=\overline{1,m}\\
c_{jj}&=\left\langle \widetilde{C}_{22},e_{j-m}^{n-m}\right\rangle ,  &c_{j,n+j}&=\left\langle \widetilde{C}_{24},e_{j-m}^{n-m}\right\rangle ,  &j &=\overline{m+1,n}
\end{aligned}\label{3.11}
\end{equation}
where $e_{j}^{m}$ denotes a $m$-dimensional canonical basis vector with $j^{th}$ component one and others zero, and $e_{j-m}^{n-m}$ is defined in the same way.

Finally, we conclude the previous deduction as 
\begin{prop}
Suppose that $A$, $B$ take the form (\ref{3.6}), (\ref{3.7}). Then the geodesics of Hamiltonian system (\ref{3.1}) that solves boundary problem (\ref{3.2}) with $t\neq\frac{k \pi}{2 a_{j}b_{j}}, \hspace{5pt} k\in \N^{+}, \hspace{5pt} j\in\overline{m+1,n}$  are given by
\begin{equation}
x(s)=
\begin{bmatrix}
C_{11} & 0 & C_{13} & 0\\
0 & C_{22} & 0 & C_{24}
\end{bmatrix}
\begin{bmatrix}
x_{1}(s)\\
x_{2}(s)\\
x_{3}(s)\\
x_{4}(s)
\end{bmatrix}
=
\begin{bmatrix}
C_{11}x_{1}(s)+C_{13}x_{3}(s)\\
C_{22}x_{2}(s)+C_{24}x_{4}(s)
\end{bmatrix}
\end{equation}
\end{prop}
where $C_{ij}$'s and components therein are identified by (\ref{3.9})-(\ref{3.11}).

\medskip
\subsection{Energy}
By use of Hamilton-Jacobi theory, the energy is conserved along the geodesics. In order to compute such energy and consequent action that both are associated with $A$, we introduce $M-$inner product $\left\langle  \cdotp, \cdotp\right\rangle _{M}:=\left\langle  M\cdotp, \cdotp\right\rangle$ with $M$ a symmetric positive definite matrix. Indeed,
\begin{align*}
&\quad \frac{d}{ds}\left[
\left\langle  \dot{x}(s),\dot{x}(s)\right\rangle _{(\varLambda^{a})^{-1}}-\left\langle  \ddot{x}(s),x(s)\right\rangle _{(\varLambda^{a})^{-1}}\right]\\
&=2\left\langle \ddot{x}(s),\dot{x}(s)\right\rangle_{(\varLambda^{a})^{-1}} -\left\langle\dddot{x}(s), x(s)\right\rangle_{(\varLambda^{a})^{-1}} -\left\langle\ddot{x}(s),\dot{x}(s) \right\rangle_{(\varLambda^{a})^{-1}} \\
&=\left\langle \ddot{x}(s),\dot{x}(s)\right\rangle_{(\varLambda^{a})^{-1}}-\left\langle\dddot{x}(s), x(s)\right\rangle_{(\varLambda^{a})^{-1}}\\
&=\left\langle Dx(s),\dot{x}(s)\right\rangle_{(\varLambda^{a})^{-1}}-\left\langle D\dot{x}(s),x(s) \right\rangle_{(\varLambda^{a})^{-1}}\\
&=0.
\end{align*}
$$
\left\langle  \dot{x}(s),\dot{x}(s)\right\rangle _{(\varLambda^{a})^{-1}}-\left\langle  \ddot{x}(s),x(s)\right\rangle _{(\varLambda^{a})^{-1}}=\mbox{Const} =:2E
$$
The main task of this subsection is to find such constant $E$ in terms of boundary data. In the following deduction, $f(T)$ denotes spectral calculus of continuous function $f$ on the selfadjoint operator $T$. As we know in the previous subsection,
$$
x(s)=
\begin{bmatrix}
C_{11} & 0 & C_{13} & 0\\
0 & C_{22} & 0 & C_{24}
\end{bmatrix}
\begin{bmatrix}
x_{1}(s)\\
x_{2}(s)\\
x_{3}(s)\\
x_{4}(s)
\end{bmatrix},
$$
\begin{align*}
\dot{x}(s)&=
\begin{bmatrix}
C_{11} & 0 & C_{13} & 0\\
0 & C_{22} & 0 & C_{24}
\end{bmatrix}\\
&\quad \cdot \begin{bmatrix}
2\sqrt{\varLambda_{1}^{a}\varLambda_{1}^{b}}& & & \\
 & & & -2\sqrt{-\varLambda_{2}^{a}\varLambda_{2}^{b}}\\
 & & -2\sqrt{\varLambda_{1}^{a}\varLambda_{1}^{b}} & \\
 & 2\sqrt{-\varLambda_{2}^{a}\varLambda_{2}^{b}} & &
\end{bmatrix}
\begin{bmatrix}
x_{1}(s)\\
x_{2}(s)\\
x_{3}(s)\\
x_{4}(s)
\end{bmatrix},
\end{align*}
and
\begin{align*}
\ddot{x}(s)=4
\begin{bmatrix}
C_{11} & 0 & C_{13} & 0\\
0 & C_{22} & 0 & C_{24}
\end{bmatrix}
\begin{bmatrix}
\varLambda_{1}^{a}\varLambda_{1}^{b} & & & \\
& \varLambda_{2}^{a}\varLambda_{2}^{b} & & \\
& & \varLambda_{1}^{a}\varLambda_{1}^{b}  & \\
& & & \varLambda_{2}^{a}\varLambda_{2}^{b}\\
\end{bmatrix}
\begin{bmatrix}
x_{1}(s)\\
x_{2}(s)\\
x_{3}(s)\\
x_{4}(s)
\end{bmatrix}.
\end{align*}

A direct computation shows
\begin{align*}
&\quad \left\langle  \dot{x}(s),\dot{x}(s)\right\rangle _{(\varLambda^{a})^{-1}}\\
&=4
\begin{bmatrix}
x_{1}(s)^{t}, & x_{2}(s)^{t}, & x_{3}(s)^{t}, & x_{4}(s)^{t}
\end{bmatrix}\\
&\quad \cdot \begin{bmatrix}
C_{11}^{2}\varLambda_{1}^{b} & 0 & -C_{11}C_{13}\varLambda_{1}^{b} & 0\\
0 & -C_{24}^{2}\varLambda_{2}^{b} & 0 & C_{22}C_{24}\varLambda_{2}^{b}\\
-C_{11}C_{13}\varLambda_{1}^{b} & 0 & C_{13}^{2}\varLambda_{1}^{b} & 0\\
0 & C_{22}C_{24}\varLambda_{2}^{b} & 0 & -C_{22}^{2}\varLambda_{2}^{b}
\end{bmatrix}
\begin{bmatrix}
x_{1}(s)\\
x_{2}(s)\\
x_{3}(s)\\
x_{4}(s)
\end{bmatrix}
\end{align*}
and 
\begin{align*}
&\quad \left\langle  \ddot{x}(s),x(s)\right\rangle _{(\varLambda^{a})^{-1}}\\
&=4
\begin{bmatrix}
x_{1}(s)^{t}, & x_{2}(s)^{t}, & x_{3}(s)^{t}, & x_{4}(s)^{t}
\end{bmatrix}\\
&\quad \cdot \begin{bmatrix}
C_{11}^{2}\varLambda_{1}^{b} & 0 & C_{11}C_{13}\varLambda_{1}^{b} & 0\\
0 & C_{22}^{2}\varLambda_{2}^{b} & 0 & C_{22}C_{24}\varLambda_{2}^{b}\\
C_{11}C_{13}\varLambda_{1}^{b} & 0 & C_{13}^{2}\varLambda_{1}^{b} & 0\\
0 & C_{22}C_{24}\varLambda_{2}^{b} & 0 & C_{24}^{2}\varLambda_{2}^{b}
\end{bmatrix}
\begin{bmatrix}
x_{1}(s)\\
x_{2}(s)\\
x_{3}(s)\\
x_{4}(s)
\end{bmatrix}.
\end{align*}
So,
\begin{align*}
&\quad
\left\langle  \dot{x}(s),\dot{x}(s)\right\rangle _{(\varLambda^{a})^{-1}}-\left\langle  \ddot{x}(s),x(s)\right\rangle _{(\varLambda^{a})^{-1}}\\
&=4
\begin{bmatrix}
x_{1}(s)^{t}, & x_{2}(s)^{t}, & x_{3}(s)^{t}, & x_{4}(s)^{t}
\end{bmatrix}\\
&\cdot\begin{bmatrix}
0 & 0 & -2C_{11}C_{13}\varLambda_{1}^{b} & 0\\
0 & -(C_{22}^{2}+C_{24}^{2})\varLambda_{2}^{b} & 0 & 0\\
-2C_{11}C_{13}\varLambda_{1}^{b} & 0 & 0 & 0\\
0 & 0 & 0 & -(C_{22}^{2}+C_{24}^{2})\varLambda_{2}^{b}
\end{bmatrix}
\begin{bmatrix}
x_{1}(s)\\
x_{2}(s)\\
x_{3}(s)\\
x_{4}(s)
\end{bmatrix}\\
&=-16 x_{3}(s)^{t}C_{11}C_{13}\varLambda_{1}^{b}x_{1}(s)\\
&\quad-4x_{2}(s)^{t}(C_{22}^{2}+C_{24}^{2})\varLambda_{2}^{b}x_{2}(s)-4x_{4}(s)^{t}(C_{22}^{2}+C_{24}^{2})\varLambda_{2}^{b}x_{4}(s)\\
&=4\{-4\tr(C_{11}C_{13}\varLambda_{1}^{b})-\tr[(C_{22}^{2}+C_{24}^{2})\varLambda_{2}^{b}]\}.
\end{align*}
Hence,
\begin{align}
E &=\frac{1}{2}(\left\langle  \dot{x},\dot{x}\right\rangle _{(\varLambda^{a})^{-1}}-\left\langle  \ddot{x},x\right\rangle _{(\varLambda^{a})^{-1}})\notag\\
&=2\left\lbrace  -4\tr(C_{11}C_{13}\varLambda_{1}^{b})-\tr\left[ (C_{22}^{2}+C_{24}^{2})\varLambda_{2}^{b}\right] \right\rbrace .
\end{align}

Making use of $\widetilde{C}_{ij}$'s solved previously, we have 
\begin{align*}
&\quad \tr\left( C_{11}C_{13}\varLambda_{1}^{b}\right) \\
&= \left\langle \widetilde{C}_{11},\widetilde{C}_{13}\right\rangle_{\varLambda_{1}^{b}}\\
&=-\left\langle  \frac{\varLambda_{1}^{b}e^{2t\sqrt{\varLambda_{1}^{a}\varLambda_{1}^{b}}}x_{1}-\varLambda_{1}^{b}x_{1}^{0}}{e^{4t\sqrt{\varLambda_{1}^{a}\varLambda_{1}^{b}}}-1}, \frac{e^{2t\sqrt{\varLambda_{1}^{a}\varLambda_{1}^{b}}}x_{1}-e^{4t\sqrt{\varLambda_{1}^{a}\varLambda_{1}^{b}}}x_{1}^{0}}{e^{4t\sqrt{\varLambda_{1}^{a}\varLambda_{1}^{b}}}-1}\right\rangle  \\
&=-\frac{1}{4}\left\langle \frac{\varLambda_{1}^{b}}{\sinh^{2}(2t\sqrt{\varLambda_{1}^{a}\varLambda_{1}^{b}})}x_{1},x_{1}\right\rangle -\frac{1}{4}\left\langle \frac{\varLambda_{1}^{b}}{\sinh^{2}(2t\sqrt{\varLambda_{1}^{a}\varLambda_{1}^{b}})}x_{1}^{0},x_{1}^{0}\right\rangle \\
&\quad+\frac{1}{2}\left\langle \frac{\varLambda_{1}^{b}\cosh(2t\sqrt{\varLambda_{1}^{a}\varLambda_{1}^{b}})}{\sinh^{2}(2t\sqrt{\varLambda_{1}^{a}\varLambda_{1}^{b}})}x_{1},x_{1}^{0}\right\rangle 
\end{align*}
and 
\begin{align*}
&\quad \tr\left[ \left( C_{24}^{2}+C_{22}^{2}\right)\varLambda_{2}^{b}\right] \\
&=\left\langle \widetilde{C}_{24},\widetilde{C}_{24}\right\rangle _{\varLambda_{2}^{b}}+\left\langle \widetilde{C}_{22},\widetilde{C}_{22}\right\rangle _{\varLambda_{2}^{b}}\\
&=\left\langle 
\varLambda_{2}^{b}\frac{x_{2}-\cos(2t\sqrt{-\varLambda_{2}^{a}\varLambda_{2}^{b}})x_{2}^{0}}{\sin(2t\sqrt{-\varLambda_{2}^{a}\varLambda_{2}^{b}})},\frac{x_{2}-\cos(2t\sqrt{-\varLambda_{2}^{a}\varLambda_{2}^{b}})x_{2}^{0}}{\sin(2t\sqrt{-\varLambda_{2}^{a}\varLambda_{2}^{b}})}\right\rangle +\left\langle \varLambda_{2}^{b}x_{2}^{0},x_{2}^{0}\right\rangle \\
&=\left\langle \frac{\varLambda_{2}^{b}}{\sin^{2}(2t\sqrt{-\varLambda_{2}^{a}\varLambda_{2}^{b}})}x_{2},x_{2}\right\rangle +\left\langle \frac{\varLambda_{2}^{b}}{\sin^{2}(2t\sqrt{-\varLambda_{2}^{a}\varLambda_{2}^{b}})}x_{2}^{0},x_{2}^{0}\right\rangle \\
&\quad-2\left\langle \frac{\varLambda_{2}^{b}\cos(2t\sqrt{-\varLambda_{2}^{a}\varLambda_{2}^{b}})}{\sin^{2}(2t\sqrt{-\varLambda_{2}^{a}\varLambda_{2}^{b}})}x_{2},x_{2}^{0}\right\rangle .
\end{align*}
Finally, we conclude the following proposition on the energy 
\begin{prop}
Suppose that $A$, $B$ take the form (\ref{3.6}), (\ref{3.7}). Then energy of Hamiltonian system (\ref{3.1}) conforms conservation law along the geodesics (\ref{3.9}) with constant $E$ given by 
\begin{equation}
\begin{split}
E&=\frac{1}{2}(\left\langle  \dot{x},\dot{x}\right\rangle _{(\varLambda^{a})^{-1}}-\left\langle  \ddot{x},x\right\rangle _{(\varLambda^{a})^{-1}})\notag\\
&=2\left\lbrace -4 \mbox{\textup{tr}}\left( C_{11}C_{13}\varLambda_{1}^{b}\right) -\mbox{\textup{tr}}\left[ \left( C_{24}^{2}+C_{22}^{2}\right) \varLambda_{2}^{b}\right] \right\rbrace \\
&=2\left[ \left\langle \frac{\varLambda_{1}^{b}}{\sinh^{2}(2t\sqrt{\varLambda_{1}^{a}\varLambda_{1}^{b}})}x_{1},x_{1}\right\rangle +\left\langle \frac{\varLambda_{1}^{b}}{\sinh^{2}(2t\sqrt{\varLambda_{1}^{a}\varLambda_{1}^{b}})}x_{1}^{0},x_{1}^{0}\right\rangle \right.\\
&\quad -2\left\langle \frac{\varLambda_{1}^{b}\cosh(2t\sqrt{\varLambda_{1}^{a}\varLambda_{1}^{b}})}{\sinh^{2}(2t\sqrt{\varLambda_{1}^{a}\varLambda_{1}^{b}})}x_{1},x_{1}^{0}\right\rangle + \left\langle \frac{-\varLambda_{2}^{b}}{\sin^{2}(2t\sqrt{-\varLambda_{2}^{a}\varLambda_{2}^{b}})}x_{2},x_{2}\right\rangle\\
&\quad+\left.\left\langle \frac{-\varLambda_{2}^{b}}{\sin^{2}(2t\sqrt{-\varLambda_{2}^{a}\varLambda_{2}^{b}})}x_{2}^{0},x_{2}^{0}\right\rangle -2\left\langle \frac{-\varLambda_{2}^{b}\cos(2t\sqrt{-\varLambda_{2}^{a}\varLambda_{2}^{b}})}{\sin^{2}(2t\sqrt{-\varLambda_{2}^{a}\varLambda_{2}^{b}})}x_{2},x_{2}^{0}\right\rangle \right]  .
\end{split}
\end{equation}
\end{prop}

\medskip
\subsection{Action function}
In this subsection, we compute the Hamilton-Jacobi action function $S$, which is a crucial ingredient in the construction of heat kernel. It satisfies \textit{Hamilton-Jacobi equation} (cf. \cite{BGG96a},~\cite{BGG96b} and \cite{B2})
$$
\frac{\partial S}{\partial t}+ H(x,\nabla S)=0.
$$
Noting in our case $H=-\frac{1}{2}E$, we have $S=\frac{1}{2}\int E dt+ c$. In the multiplier method to be adopted in section 3, the factor $\frac{1}{2}$ and constant $c$ independent of variable $t$ will be absorbed by multiplier and volume element respectively. For this reason, we do not differentiate energy from Hamiltonian, and simply define action function as 
$S=-\int E dt$. Integration by parts shows
\begin{align*}
J_{1}&= \int\frac{\varLambda_{1}^{b}}{\sinh^{2}\left( 2t\sqrt{\varLambda_{1}^{a}\varLambda_{1}^{b}}\right) } dt=-\frac{1}{2}\sqrt{\frac{\varLambda_{1}^{b}}{\varLambda_{1}^{a}}}\frac{\cosh\left( 2t\sqrt{\varLambda_{1}^{a}\varLambda_{1}^{b}}\right) }{\sinh\left( 2t\sqrt{\varLambda_{1}^{a}\varLambda_{1}^{b}}\right) },\\
J_{2}&= \int\frac{\varLambda_{1}^{b}\cosh\left( 2t\sqrt{\varLambda_{1}^{a}\varLambda_{1}^{b}}\right) }{\sinh^{2}\left( 2t\sqrt{\varLambda_{1}^{a}\varLambda_{1}^{b}}\right) }dt=-\frac{1}{2}\sqrt{\frac{\varLambda_{1}^{b}}{\varLambda_{1}^{a}}}\frac{1}{\sinh\left( 2t\sqrt{\varLambda_{1}^{a}\varLambda_{1}^{b}}\right) },\\
J_{3}&= \int\frac{-\varLambda_{2}^{b}}{\sin^{2}\left( 2t\sqrt{-\varLambda_{2}^{a}\varLambda_{2}^{b}}\right) }dt=-\frac{1}{2}\sqrt{\frac{-\varLambda_{2}^{b}}{\varLambda_{2}^{a}}}\frac{\cos\left( 2t\sqrt{-\varLambda_{2}^{a}\varLambda_{2}^{b}}\right) }{\sin\left( 2t\sqrt{-\varLambda_{2}^{a}\varLambda_{2}^{b}}\right) },\\
J_{4}&= \int\frac{-\varLambda_{2}^{b}\cos\left( 2t\sqrt{-\varLambda_{2}^{a}\varLambda_{2}^{b}}\right) }{\sin^{2}\left( 2t\sqrt{-\varLambda_{2}^{a}\varLambda_{2}^{b}}\right) }dt=-\frac{1}{2}\sqrt{\frac{-\varLambda_{2}^{b}}{\varLambda_{2}^{a}}}\frac{1}{\sin\left( 2t\sqrt{-\varLambda_{2}^{a}\varLambda_{2}^{b}}\right) }.
\end{align*}
Finally, we have 
\begin{prop}
Suppose that $A$, $B$ take the form (\ref{3.6}), (\ref{3.7}). Then the Hamilton-Jacobi action function of Hamiltonian system (\ref{3.1}) with boundary condition in (\ref{3.2}) is given by
\begin{equation}
\begin{split}
S&=-\int E dt\\
&=\left\langle  \sqrt{\frac{\varLambda_{1}^{b}}{\varLambda_{1}^{a}}}\frac{\cosh\left( 2t\sqrt{\varLambda_{1}^{a}\varLambda_{1}^{b}}\right) }{\sinh\left( 2t\sqrt{\varLambda_{1}^{a}\varLambda_{1}^{b}}\right) }x_{1},x_{1}\right\rangle +\left\langle  \sqrt{\frac{\varLambda_{1}^{b}}{\varLambda_{1}^{a}}}\frac{\cosh\left( 2t\sqrt{\varLambda_{1}^{a}\varLambda_{1}^{b}}\right) }{\sinh\left( 2t\sqrt{\varLambda_{1}^{a}\varLambda_{1}^{b}}\right) }x_{1}^{0},x_{1}^{0}\right\rangle \\
&-2\left\langle \sqrt{\frac{\varLambda_{1}^{b}}{\varLambda_{1}^{a}}}\frac{1}{\sinh\left( 2t\sqrt{\varLambda_{1}^{a}\varLambda_{1}^{b}}\right) }x_{1},x_{1}^{0}\right\rangle +\left\langle \sqrt{\frac{-\varLambda_{2}^{b}}{\varLambda_{2}^{a}}}\frac{\cos\left( 2t\sqrt{-\varLambda_{2}^{a}\varLambda_{2}^{b}}\right) }{\sin\left( 2t\sqrt{-\varLambda_{2}^{a}\varLambda_{2}^{b}}\right) }x_{2},x_{2}\right\rangle \\
&+\left\langle \sqrt{\frac{-\varLambda_{2}^{b}}{\varLambda_{2}^{a}}}\frac{\cos\left( 2t\sqrt{-\varLambda_{2}^{a}\varLambda_{2}^{b}}\right) }{\sin\left( 2t\sqrt{-\varLambda_{2}^{a}\varLambda_{2}^{b}}\right) }x_{2}^{0},x_{2}^{0}\right\rangle -2\left\langle \sqrt{\frac{-\varLambda_{2}^{b}}{\varLambda_{2}^{a}}}\frac{1}{\sin\left( 2t\sqrt{-\varLambda_{2}^{a}\varLambda_{2}^{b}}\right) }x_{2},x_{2}^{0}\right\rangle .
\end{split}\label{action}
\end{equation}
\end{prop}

\bigskip
\section{Heat kernel for $L_{S}$}
Given diagonal coefficient matrices, we find heat kernel for $L_{S}$ via multiplier techniques, and discuss its properties similar to the normal heat distribution. Heat kernel formulae for non-diagonal coefficient matrices will be given at the end of this section.

\medskip
\subsection{Explicit formulae (diagonal case)}
We start with a basic fact on the action function.

\begin{lem}
Given $A$, $B$, energy $E$ and action function $S$ as in (\ref{3.6}), (\ref{3.7}), (\ref{3.10}) and (\ref{3.11}) respectively, the following equalities hold
\begin{align}
\hspace{-10pt} & (1) \hspace{10pt} \lvert\nabla_{x} S\rvert_{A}^{2}=4\langle Bx,x\rangle + 2E\label{4.1}\\
\hspace{-10pt} & (2) \hspace{10pt} \mbox{\textup{tr}}(A \mbox{\textup{Hess}}(S))= \sum_{j=1}^{m}\frac{2a_{j}b_{j}\cosh(2t a_{j}b_{j})}{\sinh(2t a_{j}b_{j})}+ \sum_{j=m+1}^{n}\frac{2a_{j}b_{j}\cos(2t a_{j}b_{j})}{\sin(2t a_{j}b_{j})}\label{4.2}
\end{align}
where $\lvert\cdotp\rvert_{A}:=\sqrt{\langle\cdotp,\cdotp\rangle_{A}}$, and $\mbox{\textup{Hess}}(f)$ denotes Hessian of function $f\in C^{2}$.
\end{lem}

\begin{proof}
A direct computation shows that\\
(1) 
\begin{align*}
\nabla_{x} S=
\begin{bmatrix}
\nabla_{x_{1}} S\\
\nabla_{x_{2}} S
\end{bmatrix}
=2
\begin{bmatrix}
\sqrt{\frac{\varLambda_{1}^{b}}{\varLambda_{1}^{a}}}\frac{\cosh(2t\sqrt{\varLambda_{1}^{a}\varLambda_{1}^{b}})}{\sinh(2t\sqrt{\varLambda_{1}^{a}\varLambda_{1}^{b}})}x_{1}-\sqrt{\frac{\varLambda_{1}^{b}}{\varLambda_{1}^{a}}}\frac{1}{\sinh(2t\sqrt{\varLambda_{1}^{a}\varLambda_{1}^{b}})}x_{1}^{0}\\
\sqrt{\frac{-\varLambda_{2}^{b}}{\varLambda_{2}^{a}}}\frac{\cos(2t\sqrt{-\varLambda_{2}^{a}\varLambda_{2}^{b}})}{\sin(2t\sqrt{-\varLambda_{2}^{a}\varLambda_{2}^{b}})}x_{2}-\sqrt{\frac{-\varLambda_{2}^{b}}{\varLambda_{2}^{a}}}\frac{1}{\sin(2t\sqrt{-\varLambda_{2}^{a}\varLambda_{2}^{b}})}x_{2}^{0}
\end{bmatrix}
\end{align*}
\begin{align*}
&\quad\lvert\nabla_{x} S\rvert_{A}^{2}\\
&=
\begin{bmatrix}
(\nabla_{x_{1}} S)^{t},& (\nabla_{x_{2}} S)^{t}
\end{bmatrix}
\begin{bmatrix}
\varLambda_{1}^{a} & \\
 & \varLambda_{2}^{a}
\end{bmatrix}
\begin{bmatrix}
\nabla_{x_{1}} S\\
\nabla_{x_{2}} S
\end{bmatrix}\\
&=\left\langle \nabla_{x_{1}} S,\nabla_{x_{1}} S\right\rangle _{\varLambda_{1}^{a}}+\left\langle \nabla_{x_{2}} S,\nabla_{x_{2}} S\right\rangle _{\varLambda_{2}^{a}}\\
&=4\left\lbrace \left\langle \frac{\varLambda_{1}^{b}}{\varLambda_{1}^{a}}\frac{\cosh^{2}(2t\sqrt{\varLambda_{1}^{a}\varLambda_{1}^{b}})}{\sinh^{2}(2t\sqrt{\varLambda_{1}^{a}\varLambda_{1}^{b}})}x_{1},x_{1}\right\rangle _{\varLambda_{1}^{a}}+\left\langle \frac{\varLambda_{1}^{b}}{\varLambda_{1}^{a}}\frac{1}{\sinh^{2}(2t\sqrt{\varLambda_{1}^{a}\varLambda_{1}^{b}})}x_{1}^{0},x_{1}^{0}\right\rangle _{\varLambda_{1}^{a}}\right.\\
&-2\left\langle \frac{\varLambda_{1}^{b}}{\varLambda_{1}^{a}}\frac{\cosh(2t\sqrt{\varLambda_{1}^{a}\varLambda_{1}^{b}})}{\sinh^{2}(2t\sqrt{\varLambda_{1}^{a}\varLambda_{1}^{b}})}x_{1},x_{1}^{0}\right\rangle _{\varLambda_{1}^{a}}+\left\langle \frac{-\varLambda_{2}^{b}}{\varLambda_{2}^{a}}\frac{\cos^{2}(2t\sqrt{-\varLambda_{2}^{a}\varLambda_{2}^{b}})}{\sin^{2}(2t\sqrt{-\varLambda_{2}^{a}\varLambda_{2}^{b}})}x_{2},x_{2}\right\rangle _{\varLambda_{2}^{a}}\\
&+\left.\left\langle \frac{-\varLambda_{2}^{b}}{\varLambda_{2}^{a}}\frac{1}{\sin^{2}(2t\sqrt{-\varLambda_{2}^{a}\varLambda_{2}^{b}})}x_{2}^{0},x_{2}^{0}\right\rangle _{\varLambda_{2}^{a}}-2\left\langle \frac{-\varLambda_{2}^{b}}{\varLambda_{2}^{a}}\frac{\cos(2t\sqrt{-\varLambda_{2}^{a}\varLambda_{2}^{b}})}{\sin^{2}(2t\sqrt{-\varLambda_{2}^{a}\varLambda_{2}^{b}})}x_{2},x_{2}^{0}\right\rangle _{\varLambda_{2}^{a}}\right\rbrace \\
%&=4\{\langle\frac{\cosh^{2}(2t\sqrt{\varLambda_{1}^{a}\varLambda_{1}^{b}})-1}{\sinh^{2}(2t\sqrt{\varLambd%a_{1}^{a}\varLambda_{1}^{b}})}x_{1},x_{1}\rangle_{\varLambda_{1}^{b}}+\langle\frac{1-\cos^{2}(2t\sqrt{-\ %arLambda_{2}^{a}\varLambda_{2}^{b}})}{\sin^{2}(2t\sqrt{-\varLambda_{2}^{a}\varLambda_{2}^{b}})}x_{2},x_{2%}\rangle_{\varLambda_{2}^{b}}+\frac{1}{2}E\}\\
&=4\left\langle  Bx,x\right\rangle +2E.
\end{align*}
(2)
$$
\Hess(S)=2
\begin{bmatrix}
\sqrt{\frac{\varLambda_{1}^{b}}{\varLambda_{1}^{a}}}\frac{\cosh(2t\sqrt{\varLambda_{1}^{a}\varLambda_{1}^{b}})}{\sinh(2t\sqrt{\varLambda_{1}^{a}\varLambda_{1}^{b}})} & \\
 &\sqrt{\frac{-\varLambda_{2}^{b}}{\varLambda_{2}^{a}}}\frac{\cos(2t\sqrt{-\varLambda_{2}^{a}\varLambda_{2}^{b}})}{\sin(2t\sqrt{-\varLambda_{2}^{a}\varLambda_{2}^{b}})}
\end{bmatrix}.
$$
Hence, 
$$
A\Hess(S)=2
\begin{bmatrix}
\sqrt{\varLambda_{1}^{a}\varLambda_{1}^{b}}\frac{\cosh(2t\sqrt{\varLambda_{1}^{a}\varLambda_{1}^{b}})}{\sinh(2t\sqrt{\varLambda_{1}^{a}\varLambda_{1}^{b}})} & \\
 &\sqrt{-\varLambda_{2}^{a}\varLambda_{2}^{b}}\frac{\cos(2t\sqrt{-\varLambda_{2}^{a}\varLambda_{2}^{b}})}{\sin(2t\sqrt{-\varLambda_{2}^{a}\varLambda_{2}^{b}})}
\end{bmatrix}.
$$
Thus,
$$
\tr(A \Hess(S))= \sum_{j=1}^{m}\frac{2a_{j}b_{j}\cosh(2t a_{j}b_{j})}{\sinh(2t a_{j}b_{j})}+ \sum_{j=m+1}^{n}\frac{2a_{j}b_{j}\cos(2t a_{j}b_{j})}{\sin(2t a_{j}b_{j})}.
$$
\end{proof}

We expect to find the heat kernel of $ L_{S} $ in the following form
$$
K(x,x^{0};t)=V(t)e^{\kappa S(x,x^{0};t)} 
$$
where the multiplier $ \kappa $ is a real number. Making use of (\ref{4.1}) in Lemma 4.1 and noticing that
\begin{align*}
PK:&=\left( \partial_{t}+L\right) K\\
&=\left( \partial_{t}-\div\left( A\nabla\right) +\left\langle  Bx,x\right\rangle \right) \left( V(t)e^{\kappa S(x,x_{0};t)}\right) \\
&=K\left( \frac{V'}{V}-\kappa E-\kappa^{2}\lvert\nabla_{x} S\rvert_{A}^{2}-\kappa \tr(A \Hess(S))+\left\langle  Bx,x\right\rangle \right) \\
&=K\left( \frac{V'}{V}-\kappa E-4\kappa^{2}\left\langle  Bx,x\right\rangle -2\kappa^{2} E-\kappa \tr(A \Hess(S))+\left\langle  Bx,x\right\rangle \right) \\
&=0
\end{align*}
for $t>0$, we choose $\kappa=-\frac{1}{2}$ and let volume element $V(t)$  satisfy \textit{transport equation}
\begin{equation}
\frac{V(t)'}{V(t)}=\kappa \tr\left( A \Hess(S)\right).\label{4.3}
\end{equation}
Readers may consult \cite{BGG96a}, \cite{BGG96b} and \cite{B2} for more types of transport equations. By (\ref{4.2}) in Lemma 4.1, we integrate equation (\ref{4.3}) to have
\begin{equation}
V(t)=C \prod_{j=1}^{m}\left(\frac{1}{\sinh\left( 2t a_{j}b_{j}\right) }\right)^{\frac{1}{2}} \prod_{j=m+1}^{n}\left( \frac{1}{\sin\left( 2t a_{j}b_{j}\right) }\right)^{\frac{1}{2}} .
\end{equation}

The constant $C$ is determined to normalise the integral in $x-$variable of the heat kernel $K$. However, it is easy to see from (\ref{action}) that the integral is divergent if $x_{0}\neq 0$. $K(x,0;t)$, a propagator from origin to arbitrary point $x$, is called \textit{generalised heat kernel} as we shall prove next section that it indeed has similar properties to the normal heat distribution. We denote $K(x,0;t)$ by $K(x;t)$ from now on, then
\begin{align*}
K(x;t)&=C\prod_{j=1}^{m}\left( \frac{1}{\sinh(2t a_{j}b_{j})}\right) ^{\frac{1}{2}}\prod_{j=m+1}^{n}\left( \frac{1}{\sin(2t a_{j}b_{j})}\right)^{\frac{1}{2}} \\
&\quad\times e^{-\frac{1}{4t}\left( \sum_{j=1}^{m}\frac{2t b_{j}}{a_{j}}\frac{\cosh(2t a_{j}b_{j})}{\sinh(2t a_{j}b_{j})}(x_{j}^{(1)})^{2}+\sum_{j=m+1}^{n}\frac{2t b_{j}}{a_{j}}\frac{\cos(2t a_{j}b_{j})}{\sin(2t a_{j}b_{j})}(x_{j}^{(1)})^{2}\right) }.
\end{align*}

Making use of $\int_{\R}e^{-x^{2}}dx=\sqrt{\pi}$,
\begin{align*}
&\quad \int_{\R^{n}}K(x;t)dx\\ 
&= C\prod_{j=1}^{m}\left( \frac{1}{\sinh(2t a_{j}b_{j})}\right) ^{\frac{1}{2}}\prod_{j=m+1}^{n}\left( \frac{1}{\sin(2t a_{j}b_{j})}\right)^{\frac{1}{2}} \\
&\quad\cdot\prod_{j=1}^{m}\left[ \left( \frac{1}{4t}\frac{2t b_{j}}{a_{j}}\frac{\cosh(2t a_{j}b_{j})}{\sinh(2t a_{j}b_{j})}\right) ^{-\frac{1}{2}}\sqrt{\pi}\right]
\prod_{j=m+1}^{n}\left[ \left( \frac{1}{4t}\frac{2t b_{j}}{a_{j}}\frac{\cos(2t a_{j}b_{j})}{\sin(2t a_{j}b_{j})}\right) ^{-\frac{1}{2}}\sqrt{\pi}\right]\\
&=C \prod_{j=1}^{n}\left( \dfrac{2\pi a_{j}}{b_{j}}\right) ^{\frac{1}{2}}\prod_{j=1}^{m}\left(\frac{1}{\cosh(2t a_{j}b_{j})} \right)^{\frac{1}{2}}\prod_{j=m+1}^{n}\left(\frac{1}{\cos(2t a_{j}b_{j})} \right)^{\frac{1}{2}}
\end{align*}
tends to $C\prod_{j=1}^{n}\left( \dfrac{2\pi a_{j}}{b_{j}}\right) ^{\frac{1}{2}}$ as $t\rightarrow 0^{+}$.

By choosing $C=\prod_{j=1}^{n}\left( \dfrac{b_{j}}{2\pi a_{j}}\right) ^{\frac{1}{2}}$, we have arrived at the following proposition
\begin{prop}
Assume that $A$ and $B$ are diagonal matrices as in (\ref{3.3}) and (\ref{3.4}). Then the heat kernel of the Schr\"{o}dinger operator $L_{S}=-\mbox{\textup{div}}(A\nabla)+\left\langle  Bx,x \right\rangle $ is
\begin{equation}
\begin{split}
K(x,x^{0};t)&=(4\pi t)^{-\frac{n}{2}}\prod_{j=1}^{m}\left( \frac{2t b_{j}}{a_{j}\sinh(2t a_{j}b_{j})}\right) ^{\frac{1}{2}}\prod_{j=m+1}^{n}\left( \frac{2t b_{j}}{a_{j}\sin(2t a_{j}b_{j})}\right)^{\frac{1}{2}} \\
&\quad\times e^{-\frac{1}{4t}\left( \sum_{j=1}^{m}\frac{2t b_{j}}{a_{j}}\frac{\cosh(2t a_{j}b_{j})}{\sinh(2t a_{j}b_{j})}(x_{j}^{(1)})^{2}+\sum_{j=1}^{m}\frac{2t b_{j}}{a_{j}}\frac{\cosh(2t a_{j}b_{j})}{\sinh(2t a_{j}b_{j})}(x_{j}^{(0)})^{2}\right) }\\
&\quad\times e^{-\frac{1}{4t}\left( \sum_{j=m+1}^{n}\frac{2t b_{j}}{a_{j}}\frac{\cos(2t a_{j}b_{j})}{\sin(2t a_{j}b_{j})}(x_{j}^{(1)})^{2}+\sum_{j=m+1}^{n}\frac{2t b_{j}}{a_{j}}\frac{\cos(2t a_{j}b_{j})}{\sin(2t a_{j}b_{j})}(x_{j}^{(0)})^{2}\right) }\\
&\quad\times e^{\frac{1}{2t}\left(\sum_{j=1}^{m}\frac{2t b_{j}}{a_{j}}\frac{1}{\sinh(2t a_{j}b_{j})}x_{j}^{(0)}x_{j}^{(1)}+\sum_{j=m+1}^{n}\frac{2t b_{j}}{a_{j}}\frac{1}{\sin(2t a_{j}b_{j})}x_{j}^{(0)}x_{j}^{(1)}\right) }.
\end{split}\label{4.5}
\end{equation}
\end{prop}

\begin{rmk}
In sake of continuity of $a_{j}$'s and $b_{j}$'s, heat kernel (\ref{4.5}) keeps valid if the matrix $D$ has multiple eigenvalues or zero eigenvalues. Moreover, condition (C) is technical, and can be removed.
\end{rmk}

\medskip
\begin{rmk}
Heat kernel (\ref{4.5}) is complex valued as long as $\sin(2t a_{j}b_{j})< 0$, i.e. $t\in (\frac{4k+1}{4a_{j}b_{j}}\pi,\frac{4k+3}{4a_{j}b_{j}}\pi )$, $k\in \N^{+}$, $j\in\overline{m+1,n}$.
\end{rmk}

\medskip
\begin{rmk}
Heat kernel (\ref{4.5}) holds if the sub-matrix  $\sin\left( 2t\sqrt{-\varLambda_{2}^{a}\varLambda_{2}^{b}}\right) $ is non-singular, which we proposed as an assumption in the previous section. We call region $\Omega=\left\lbrace \right(x,t)\in \R^{n}\times\R^{+}: t=\frac{k \pi}{2 a_{j}b_{j}}, \hspace{5pt} k\in \N^{+}, \hspace{5pt} j\in\overline{m+1,n}\rbrace $ singular region and region $\Omega^{c}=\R^{n}\times\R^{+}\setminus\Omega$ regular region. Briefly speaking, there is no geodesic or uncountably many geodesics connecting the given boundary points $x$ and $x^{0}$ for $t=\frac{k \pi}{2 a_{j}b_{j}}$, while there is a unique geodesic for any given two points $x$ and $x^{0}$ if $t\neq \frac{k \pi}{2 a_{j}b_{j}}$. Here we point out that such singular region has no contribution to the Hamilton-Jacobi action function which is regarded as an integral of energy in $t-$variable.
\end{rmk}

\medskip
\subsection{Generalised heat kernel}
In this subsection, we show that \textit{generalised heat kernel} has analogue properties to the normal one, that is 

\begin{prop}
Heat kernel (\ref{4.5}) is said to be generalised in the following sense.
\begin{align*}
& (1) \hspace{10pt} K(x,x^{0};t)>0, \hspace{5pt} \forall (x,x^{0})\in \R^{n}\times\R^{n}, \hspace{5pt} 0<t\ll 1.\\
& (2) \hspace{10pt} Fix \hspace{5pt} x^{0}=0, \hspace{5pt} \hat{K}(\xi;t)\rightarrow 1, \hspace{5pt} as \hspace{5pt} t\rightarrow 0^{+}.\\
& (3) \hspace{10pt} Fix \hspace{5pt} x^{0}=0, \hspace{5pt} K(x;t)\xrightarrow{d} \delta(x), \hspace{5pt} as \hspace{5pt} t\rightarrow 0^{+}.
\end{align*}
where hat denotes Fourier transform on spatial variables, and $\xrightarrow{d}$ means limitation in the sense of distribution.
\end{prop}

\begin{proof}
(1) It is obvious from formulae (\ref{4.5}) for $t$ appropriately small.

\bigskip
In the rest of this proof, we fix $x^{0}=0$ in (\ref{4.5}) and 
\begin{equation}
\begin{aligned}
K(x;t)&=(4\pi t)^{-\frac{n}{2}}\prod_{j=1}^{m}\left( \frac{2t b_{j}}{a_{j} \sinh(2t a_{j}b_{j})}\right) ^{\frac{1}{2}}\prod_{j=m+1}^{n}\left( \frac{2t b_{j}}{ a_{j} \sin(2t a_{j}b_{j})}\right)^{\frac{1}{2}} \\
&\quad\times e^{-\frac{1}{4t}\left( \sum_{j=1}^{m}\frac{2t b_{j}}{a_{j}}\frac{\cosh(2t a_{j}b_{j})}{\sinh(2t a_{j}b_{j})}(x_{j}^{(1)})^{2}+\sum_{j=m+1}^{n}\frac{2t b_{j}}{a_{j}}\frac{\cos(2t a_{j}b_{j})}{\sin(2t a_{j}b_{j})}(x_{j}^{(1)})^{2}\right) }.
\end{aligned}\label{4.5.1}
\end{equation}

\medskip
\noindent 
(2) By properties of Fourier transform 
$$
\widehat{e^{-\pi  x^{2}}}(\xi)=e^{-\pi \xi^{2}}
$$
and
$$
\widehat{f(\lambda x)}(\xi)=\lambda^{-1}\hat{f}\left( \lambda^{-1}\xi\right),
$$
we have
\begin{align*}
\hat{K}(\xi;t)&=\int_{\R^{n}}K(x;t)e^{-2\pi i \xi \cdot x}dx\\
&= (4\pi t)^{-\frac{n}{2}}\prod_{j=1}^{m}\left( \frac{2t b_{j}}{a_{j} \sinh(2t a_{j}b_{j})}\right) ^{\frac{1}{2}}\prod_{j=m+1}^{n}\left( \frac{2t b_{j}}{ a_{j} \sin(2t a_{j}b_{j})}\right)^{\frac{1}{2}} \\
&\quad\times\prod_{j=1}^{m}\int_{\R}e^{-\frac{1}{4t}\frac{2t b_{j}}{a_{j}}\frac{\cosh(2t a_{j}b_{j})}{\sinh(2t a_{j}b_{j})}(x_{j}^{(1)})^{2}}e^{-2\pi i \xi_{j} \cdot x_{j}^{(1)}}dx_{j}^{(1)}\\
&\quad\times \prod_{j=m+1}^{n}\int_{\R}e^{-\frac{1}{4t}\frac{2t b_{j}}{a_{j}}\frac{\cos(2t a_{j}b_{j})}{\sin(2t a_{j}b_{j})}(x_{j}^{(1)})^{2}}e^{-2\pi i \xi_{j} \cdot x_{j}^{(1)}}dx_{j}^{(1)}\\
&= (4\pi t)^{-\frac{n}{2}}\prod_{j=1}^{m}\left( \frac{2t b_{j}}{a_{j} \sinh(2t a_{j}b_{j})}\right) ^{\frac{1}{2}}\prod_{j=m+1}^{n}\left( \frac{2t b_{j}}{ a_{j} \sin(2t a_{j}b_{j})}\right)^{\frac{1}{2}} \\
&\quad\times\prod_{j=1}^{m}\left(\frac{1}{4\pi t}\frac{2t b_{j}}{a_{j}}\frac{\cosh(2t a_{j}b_{j})}{\sinh(2t a_{j}b_{j})}\right)^{-\frac{1}{2}}e^{-4\pi^{2}t\frac{a_{j}}{2t b_{j}}\frac{\sinh(2t a_{j}b_{j})}{\cosh(2t a_{j}b_{j})}\xi_{j}^{2}}\\
&\quad\times\prod_{j=m+1}^{n}\left(\frac{1}{4\pi t}\frac{2t b_{j}}{a_{j}}\frac{\cos(2t a_{j}b_{j})}{\sin(2t a_{j}b_{j})}\right)^{-\frac{1}{2}}e^{-4\pi^{2}t\frac{a_{j}}{2t b_{j}}\frac{\sin(2t a_{j}b_{j})}{\cos(2t a_{j}b_{j})}\xi_{j}^{2}}\\
&=\prod_{j=1}^{m}\left(\frac{1}{\cosh(2t a_{j}b_{j})} \right)^{\frac{1}{2}}\prod_{j=m+1}^{n}\left(\frac{1}{\cos(2t a_{j}b_{j})} \right)^{\frac{1}{2}}\\
&\quad\times\prod_{j=1}^{m}e^{-2\pi^{2}\frac{a_{j}}{b_{j}}\frac{\sinh(2t a_{j}b_{j})}{\cosh(2t a_{j}b_{j})}\xi_{j}^{2}}\prod_{j=m+1}^{n}e^{-2\pi^{2}\frac{a_{j}}{ b_{j}}\frac{\sin(2t a_{j}b_{j})}{\cos(2t a_{j}b_{j})}\xi_{j}^{2}}\rightarrow 1
\end{align*}
as $t\rightarrow 0^{+}$.

\medskip
\noindent 
(3) We write (\ref{4.5.1}) as $K(x;t)=K_{1}(x;t)K_{2}(x;t)$, where
\begin{align*}
& \qquad K_{1}(x;t)=\prod_{j=1}^{m}\left( \frac{2t a_{j}b_{j}}{\sinh(2t a_{j}b_{j})}\right) ^{\frac{1}{2}}\prod_{j=m+1}^{n}\left( \frac{2t a_{j}b_{j}}{\sin(2t a_{j}b_{j})}\right)^{\frac{1}{2}} \\ 
& \cdot e^{-\frac{1}{4t}\left[ \sum_{j=1}^{m}\left( \frac{x_{j}^{(1)}}{a_{j}}\right) ^{2}\left( \frac{2t a_{j}b_{j}\cosh(2t a_{j}b_{j})}{\sinh(2t a_{j}b_{j})}-1 \right) +\sum_{j=m+1}^{n}\left( \frac{x_{j}^{(1)}}{a_{j}}\right) ^{2}\left( \frac{2t a_{j}b_{j}\cos(2t a_{j}b_{j})}{\sin(2t a_{j}b_{j})}-1 \right)\right] }
\end{align*}
and
\begin{equation*}
K_{2}(x;t)=\frac{\left(4 \pi t \right)^{-\frac{n}{2}} }{\left( \det A \right)^{\frac{1}{2}} }e^{-\frac{\lvert x\rvert_{A^{-1}}^{2}}{4t}}.
\end{equation*}

\bigskip
The proof will be carried out in two steps.

\medskip
(i) $K_{2}(x;t) \xrightarrow{d} \delta (x)$, as $t\rightarrow 0^{+}$.

\medskip
For any $\varphi \in C_{0}^{\infty}(\R^{n})$,
\begin{equation}
\int_{\R^{n}}K_{2}(x;t)\varphi(x)dx=\varphi(0)\int_{\R^{n}}K_{2}(x;t)dx+\int_{\R^{n}}K_{2}(x;t)\left[ \varphi(x)-\varphi(0)\right]dx. \label{*}
\end{equation}

The first term 
$$
\int_{\R^{n}}K_{2}(x;t)dx=\prod_{j=1}^{n}\int_{\R}\dfrac{1}{a_{j}\sqrt{4 \pi t}}e^{-\frac{x_{j}^{2}}{4t a_{j}^{2}}}d x_{j}=\prod_{j=1}^{n}\int_{\R}e^{-\pi x_{j}^{2}}d x_{j}=1,
$$
and the second term
\begin{align*}
\int_{\R^{n}}K_{2}(x;t)\left[ \varphi(x)-\varphi(0)\right]dx &= \int_{\R^{n}}\frac{\left(4 \pi t \right)^{-\frac{n}{2}} }{\left( \det A \right)^{\frac{1}{2}} }e^{-\frac{\lvert A^{-\frac{1}{2}}x\rvert^{2}}{4t}}\left[ \varphi(x)-\varphi(0)\right]dx\\
&= \int_{\R^{n}}e^{-\pi\lvert y \rvert^{2}}\left[ \varphi\left( \sqrt{4\pi t }A^{\frac{1}{2}}y\right) -\varphi(0) \right] dy \rightarrow 0,
\end{align*}
as $t\rightarrow 0^{+}$.

As a result of (\ref{*}), 
$$
\lim_{t\rightarrow 0^{+}}\int_{\R^{n}}K_{2}(x;t)\varphi(x)dx = \varphi(0) = \int_{\R^{n}}\delta(x) \varphi(x) dx,
$$
i.e. $K_{2}\xrightarrow{d} \delta $, as $t\rightarrow 0^{+}$.

\medskip
(ii) $K(x;t) \xrightarrow{d} \delta (x)$, as $t\rightarrow 0^{+}$.

\medskip
Taking any $\varphi \in C_{0}^{\infty}(\R^{n})$, we assume that support of $\varphi$ is contained in the ball $\left\lbrace x\in \R^{n}:\, \rvert x \rvert \leqslant R\right\rbrace $, and that $\varphi$ is dominated by some constant $C$ everywhere. For the sake of 
$$
\int_{\R^{n}}K\varphi dx -\varphi(0) = \int_{\R^{n}}(K_{1}-1)K_{2}\varphi dx + \int_{\R^{n}}K_{2}\varphi dx - \varphi(0)
$$
and step (i) $\lim_{t\rightarrow 0^{+}}\int_{\R^{n}}K_{2}\varphi dx = \varphi(0)$, it is sufficient to conclude $K \xrightarrow{d} \delta$, as $t\rightarrow 0^{+}$ by checking 
$$
\lim_{t\rightarrow 0^{+}}\int_{\R^{n}}(K_{1}-1)K_{2}\varphi dx = 0.
$$
Indeed, a variable change $y=\frac{A^{-\frac{1}{2}}x}{\sqrt{4\pi t}}$ makes 
\begin{align*}
\int_{\R^{n}}(K_{1}-1)K_{2}\varphi dx &= \int_{\rvert x \rvert \leqslant R}(K_{1}(x;t)-1)K_{2}(x;t)\varphi (x) dx\\
&= \int_{\rvert y \rvert \leqslant \frac{\lVert A^{-\frac{1}{2}} \rVert R}{\sqrt{4\pi t}}}(K_{1}(\sqrt{4\pi t}A^{\frac{1}{2}}y)-1)\varphi(\sqrt{4\pi t}A^{\frac{1}{2}}y) e^{-\pi \lvert y \rvert^{2}} dy
\end{align*}
where 
\begin{align*}
& \qquad K_{1}(\sqrt{4\pi t}A^{\frac{1}{2}}y) = \prod_{j=1}^{m}\left( \frac{2t a_{j}b_{j}}{\sinh(2t a_{j}b_{j})}\right) ^{\frac{1}{2}}\prod_{j=m+1}^{n}\left( \frac{2t a_{j}b_{j}}{\sin(2t a_{j}b_{j})}\right)^{\frac{1}{2}} \\ 
& \cdot e^{-\pi \left[ \sum_{j=1}^{m} y_{j}^{2}\left( \frac{2t a_{j}b_{j}\cosh(2t a_{j}b_{j})}{\sinh(2t a_{j}b_{j})}-1 \right) +\sum_{j=m+1}^{n}y_{j}^{2} \left( \frac{2t a_{j}b_{j}\cos(2t a_{j}b_{j})}{\sin(2t a_{j}b_{j})}-1 \right)\right] }.
\end{align*}
Noticing that $\frac{u\cosh(u)-\sinh(u)}{\sinh(u)}\rightarrow 0^{+}$ and $\frac{u\cos(u)-\sin(u)}{\sin(u)}\rightarrow 0^{-}$ as $t\rightarrow 0^{+}$, we obtain 
\begin{align*}
\left|\int_{\R^{n}}(K_{1}-1)K_{2}\varphi dx \right| &\leqslant C \prod_{j=1}^{m}\left( \frac{2t a_{j}b_{j}}{\sinh(2t a_{j}b_{j})}\right) ^{\frac{1}{2}}\prod_{j=m+1}^{n}\left( \frac{2t a_{j}b_{j}}{\sin(2t a_{j}b_{j})}\right)^{\frac{1}{2}} \\ 
& \cdotp \int_{\R^{n}}\left[ e^{-\pi \sum_{j=m+1}^{n}y_{j}^{2} \left( \frac{2t a_{j}b_{j}\cos(2t a_{j}b_{j})}{\sin(2t a_{j}b_{j})}-1 \right)} -1 \right] e^{-\pi \lvert y \rvert^{2}} dy \rightarrow 0
\end{align*}
as $t\rightarrow 0^{+}$, which completes the proof.
\end{proof}

\medskip
\subsection{Explicit formulae (non-diagonal case)}
In order to generalise Proposition 3.1 to non-diagonal case, we first rewrite (\ref{4.5}) in inner-product form
\begin{equation}
\begin{split}
&\quad K(x,x^{0};t)=(4\pi t)^{-\frac{n}{2}}\left( \frac{\det\psi\left( 2t\sqrt{\varLambda^{a}\varLambda^{b}}\right)}{\det \varLambda^{a}}  \right) ^{\frac{1}{2}}\\
&\times e^{-\frac{1}{4t}\left(\left\langle \varphi\left(2t\sqrt{\varLambda^{a}\varLambda^{b}} \right) x,x \right\rangle _{(\varLambda^{a})^{-1}}+\left\langle \varphi\left(2t\sqrt{\varLambda^{a}\varLambda^{b}} \right) x^{0},x^{0} \right\rangle _{(\varLambda^{a})^{-1}}-2\left\langle \psi\left(2t\sqrt{\varLambda^{a}\varLambda^{b}} \right) x,x^{0} \right\rangle _{(\varLambda^{a})^{-1}} \right) }
\end{split}\label{4.6}
\end{equation}
where $\varphi(u)=u\coth(u)$, $\psi(u)=\frac{u}{\sinh(u)}$.

Suppose that $A$ and $B$ are commutative, then there exists an orthogonal matrix $P$ such that $PAP^{t}=\varLambda^{a}$ and $P(B+B^{t})P^{t}=2\varLambda^{b}$. Putting $y=Px$, system $\ddot{x}=2A(B+B^{t})x=:Dx$ becomes $\ddot{y}=4\varLambda^{a}\varLambda^{b}y$. Moreover, we have 
$$
\det\psi\left( 2t\sqrt{\varLambda^{a}\varLambda^{b}}\right)= \det\psi\left(t\sqrt{D}\right),
$$
\begin{align*}
\left\langle \varphi\left(2t\sqrt{\varLambda^{a}\varLambda^{b}} \right) y,y \right\rangle _{(\varLambda^{a})^{-1}}&=\left\langle(\varLambda^{a})^{-1}\varphi\left(2t\sqrt{\varLambda^{a}\varLambda^{b}} \right)Px,Px \right\rangle\\ 
&=\left\langle PA^{-1}P^{t}\varphi\left(2t\sqrt{\varLambda^{a}\varLambda^{b}} \right)Px,Px \right\rangle\\
&=\left\langle A^{-1}\varphi\left(t\sqrt{D}\right)x, x \right\rangle\\
&=\left\langle \varphi\left(t\sqrt{D}\right)x, x \right\rangle _{A^{-1}},
\end{align*}
and similarly,
$$
\left\langle \varphi\left(2t\sqrt{\varLambda^{a}\varLambda^{b}} \right) y^{0},y^{0} \right\rangle _{(\varLambda^{a})^{-1}}=\left\langle \varphi\left(t\sqrt{D}\right)x^{0}, x^{0} \right\rangle _{A^{-1}},
$$
$$
\left\langle \psi\left(2t\sqrt{\varLambda^{a}\varLambda^{b}} \right) y,y^{0} \right\rangle _{(\varLambda^{a})^{-1}}=\left\langle \psi\left(t\sqrt{D}\right)x, x^{0} \right\rangle _{A^{-1}}.
$$
According to (\ref{4.6}), we have arrived one of our main results.
\begin{thm}
For given matrix $A$ symmetric positive definite, $B$ a real matrix such that $A$ and $B$ are commutative, the heat kernel of Schr\"{o}dinger operator $$L_{S}=-\mbox{\textup{div}}(A\nabla)+\langle Bx,x \rangle$$ has the following form
\begin{align}
K(x,x^{0};t)&=(4\pi t)^{-\frac{n}{2}}\left(\frac{\det\psi\left(t\sqrt{D}\right)}{\det A} \right) ^{\frac{1}{2}}\notag\\
&\quad\times e^{-\frac{1}{4t}\left(\left\langle \varphi\left(t\sqrt{D}\right)x, x \right\rangle _{A^{-1}}+\left\langle \varphi\left(t\sqrt{D}\right)x^{0}, x^{0} \right\rangle _{A^{-1}}-2\left\langle \psi\left(t\sqrt{D}\right)x, x^{0} \right\rangle _{A^{-1}}\right)}\tag{2.1}
\end{align} 
where $D=2A(B+B^{t})$, $\varphi(u)=u\coth(u)$, $\psi(u)=\frac{u}{\sinh(u)}$.
\end{thm}

\medskip
\begin{rmk}
We use the convention $u\coth(u)\arrowvert_{u=0}=1$, so $D\coth(D)=I_{n}$ on the kernel of $D$, and so on.
\end{rmk}

\bigskip
\section{Heat kernel for $L$}
For the first application of Theorem 3.1, we compute the explicit heat kernel of operator $L$. With mention in introduction, we will adopt the ansatz (\ref{1.3}) and solve the associated differential equation system (\ref{1.4})-(\ref{1.9}). Of all these equations, the most difficult one is the matrix Riccati equation (\ref{1.4}). Given the results of previous sections, we handle this point in a straightforward way. In fact, we have the following 
\begin{thm} \textbf{\mbox{\textup{(Globally closed solution for matrix Riccati equation)}}}
For matrix $A$ symmetric positive definite, $B$ a real matrix such that $A$ and $B$ are commutative, 
matrix Riccati equation
\begin{equation}
\dot{\alpha}=4\alpha A\alpha-\frac{B+B^{t}}{2} \tag{1.4}
\end{equation}
has a globally explicit solution
\begin{equation}
\alpha =\frac{-1}{4t}A^{-1}\varphi\left(t\sqrt{D}\right).
\end{equation}
Besides, the differential equations system 
\begin{align}
& \dot{\beta}=4\beta A\alpha \tag{1.5} \\
& \dot{\gamma}=\beta A\beta  \tag{1.6}
\end{align}
have explicit solutions
\begin{equation}
\begin{split}
\beta &=\frac{1}{2t}A^{-1}\psi\left(t\sqrt{D}\right)\\
\gamma &=\frac{-1}{4t}A^{-1}\varphi\left(t\sqrt{D}\right)
\end{split}\label{4.7}
\end{equation}
where $D=2A(B+B^{t})$, $\varphi(u)=u\coth(u)$, $\psi(u)=\frac{u}{\sinh(u)}$.
\end{thm}

\begin{proof}
The solutions are read out from the Theorem 3.1 since equations (\ref{1.4})-(\ref{1.6}) are irrelevant to vector $f$, $g$ and constant $h$. 
\end{proof}

Next, we may integrate equations (\ref{1.7})-(\ref{1.9}) for $A$ and $B$ satisfying condition of Theorem 3.1. We point out that for singular matrix $B$, solutions are formulated in component form. To have concise solutions, we assume that $B$ is non-singular and commutative with symmetric positive definite $A$.

\medskip
$\bullet \hspace{10pt}\mu-function$
\begin{equation}
\mu=\frac{1}{2}A^{-1}f-\frac{\cosh\left( t\sqrt{D}\right)}{\sqrt{D}\sinh\left( t\sqrt{D}\right)}g.\label{4.8}
\end{equation}

$\bullet \hspace{10pt}\nu-function$
\begin{equation}
\nu=\frac{1}{\sqrt{D}\sinh\left( t\sqrt{D}\right)}g.
\end{equation}

\begin{rmk}
Free constants in $\mu-$function and $\nu-$function are absorbed by function $W(t)$.
\end{rmk}

$\bullet \hspace{10pt}W-function$\\
Making use of 
\begin{align*}
& A\alpha =\frac{-1}{4t}\varphi\left( t\sqrt{D}\right),\\
& \left\langle A\mu,\mu\right\rangle=\frac{1}{4}\lvert f\rvert_{A^{-1}}^{2}-\left\langle \frac{\coth\left( t\sqrt{D}\right)}{\sqrt{D}}f,g\right\rangle +\left\langle\frac{\coth^{2}\left( t\sqrt{D}\right)}{2\left(B+B^{t} \right) }g,g \right\rangle,\\
& \left\langle f, \mu \right\rangle =\frac{1}{2}\lvert f\rvert_{A^{-1}}^{2}-\left\langle \frac{\coth\left( t\sqrt{D}\right)}{\sqrt{D}}f,g\right\rangle,
\end{align*} 
volume element $W(t)$ satisfies
$$
W^{-1}\dot{W}=\frac{-1}{2t}\tr\left[\varphi \left( t\sqrt{D}\right) \right]-\frac{1}{4}\lvert f\rvert_{A^{-1}}^{2}+\left\langle\frac{\coth^{2}\left( t\sqrt{D}\right)}{2\left(B+B^{t} \right) }g,g \right\rangle-h.
$$
Integration yields
$$
W(t)=W_{0}e^{-\left( \frac{1}{4}\lvert f\rvert_{A^{-1}}^{2}+h\right)t+\left\langle\phi\left(t\sqrt{D} \right)g, g  \right\rangle_{A} t^{3}}
$$ 
where $W_{0}=C\left[ \det\frac{1}{\sinh\left( t\sqrt{D}\right)}\right] ^{\frac{1}{2}}$ and $\phi(u)=\frac{u-\coth(u)}{u^{3}}$. Taking $f=g=h=0$, we have 
$$
W_{0}=W=V=(4\pi t)^{-\frac{n}{2}}\left(\frac{\det\psi\left(t\sqrt{D}\right)}{\det A} \right) ^{\frac{1}{2}}.
$$
Consequently, volume element $W$ is given by
\begin{equation}
W=(4\pi t)^{-\frac{n}{2}}\left(\frac{\det\psi\left(t\sqrt{D}\right)}{\det A} \right) ^{\frac{1}{2}} e^{-\left( \frac{1}{4}\lvert f\rvert_{A^{-1}}^{2}+h\right)t+\left\langle\phi\left(t\sqrt{D} \right)g, g  \right\rangle_{A} t^{3}}.
\end{equation}

Finally, we yields another main result of this paper. 
\begin{thm}
For $A$ symmetric positive definite, $B$ non-singular such that $A$ and $B$ are commutative, the heat kernel of operator 
$$
L=-\mbox{\textup{div}}(A\nabla) + \langle Bx,x \rangle + \langle f,\nabla \rangle + \langle g,x \rangle+h
$$ 
has the following form
\begin{align}
K(x,x^{0};t)&=(4\pi t)^{-\frac{n}{2}}\left(\frac{\det\psi\left(t\sqrt{D}\right)}{\det A} \right) ^{\frac{1}{2}}e^{-\left( \frac{1}{4}\lvert f\rvert_{A^{-1}}^{2}+h\right)t+\left\langle\phi\left(t\sqrt{D} \right)g, g  \right\rangle_{A} t^{3}}\notag\\
&\quad\times e^{-\frac{1}{4t}\left(\left\langle \varphi\left(t\sqrt{D}\right)x, x \right\rangle _{A^{-1}}+\left\langle \varphi\left(t\sqrt{D}\right)x^{0}, x^{0} \right\rangle _{A^{-1}}-2\left\langle \psi\left(t\sqrt{D}\right)x, x^{0} \right\rangle _{A^{-1}}\right)}\notag\\
&\quad\times e^{\frac{1}{2}\left\langle f,x \right\rangle _{A^{-1}}-\langle \frac{\coth\left( t\sqrt{D}\right)}{\sqrt{D}}g, x\rangle +\langle \frac{1}{\sqrt{D}\sinh\left( t\sqrt{D}\right)}g, x^{0}\rangle }\tag{2.2}
\end{align}
where $D=2A(B+B^{t})$, $\varphi(u)=u\coth(u)$, $\psi(u)=\frac{u}{\sinh(u)}$, $\phi(u)=\frac{u-\coth(u)}{u^{3}}$.
\end{thm} 

\bigskip
\section{Examples}
For the second application of Theorem 3.1, we demonstrate three examples to recover and generalise  several classical results on some celebrated operators. Notation $\varphi(u)=u\coth(u)$ and $\psi(u)=\frac{u}{\sinh(u)}$ keep valid throughout the whole section. 

\medskip
\begin{eg}
Generalised Laplacian
\end{eg}
\noindent
Define generalised Laplacian as $L_{GL}=-\div(A\nabla)$ where $A$ is a symmetric positive definite matrix. With $B=0$ in Theorem 3.1, we have $D=0$, $\varphi\left(t\sqrt{D}\right)=I_{n}$, $\psi\left(t\sqrt{D}\right)=I_{n}$ and the heat kernel is given by
\begin{equation}
\begin{split}
K_{GL}(x,x^{0};t)&=\frac{(4\pi t)^{-\frac{n}{2}}}{\left(\det A \right)^{\frac{1}{2}} }e^{-\frac{1}{4t}\left(\left\langle x,x\right\rangle_{A^{-1}} +\left\langle x^{0}, x^{0}\right\rangle_{A^{-1}}-2\left\langle x, x^{0}\right\rangle_{A^{-1}}\right)}\\
&=\frac{(4\pi t)^{-\frac{n}{2}}}{\left(\det A \right)^{\frac{1}{2}} }e^{-\frac{\lvert x- x^{0} \rvert_{A^{-1}}^{2}}{4t}}.
\end{split}
\end{equation}
In particular, taking $A=I_{n}$, kernel becomes
$$
K_{L}(x,x^{0};t)=(4\pi t)^{-\frac{n}{2}}e^{-\frac{\lvert x- x^{0} \rvert^{2}}{4t}}
$$
which is exactly the Gaussian.

\medskip
\begin{eg}
Generalised Hermite operator
\end{eg}
\noindent
Define generalised Hermite operator $L_{GH}=-\div(A\nabla)+\left\langle  Bx,x \right\rangle $ with $B>0$.\\
\medskip
\hspace{15pt}$\bullet \hspace{10pt} A=I_{n}$, \hspace{5pt}$B=\diag \{b_{j}^{2}\}_{j=1}^{n}$  \hspace{10pt}($b_{j}>0$)\\
By Theorem 3.1, we yield Mehler formulae
\begin{equation}
\begin{split}
K(x,x^{0};t)&=(4\pi t)^{-\frac{n}{2}}\prod_{j=1}^{n}\left( \frac{2t b_{j}}{\sinh(2t b_{j})}\right) ^{\frac{1}{2}}\\
&\quad\times e^{-\frac{1}{4t}\left( \sum_{j=1}^{n}\frac{2t b_{j}\cosh(2t b_{j})}{\sinh(2t b_{j})}(x_{j}^{(1)})^{2}+\sum_{j=1}^{n}\frac{2t b_{j}\cosh(2t b_{j})}{\sinh(2t b_{j})}(x_{j}^{(0)})^{2}\right) }\\
&\quad\times e^{-\frac{1}{4t}\left(-2\sum_{j=1}^{n}\frac{2t b_{j}}{\sinh(2t b_{j})}x_{j}^{(0)}x_{j}^{(1)}\right) }.
\end{split}
\end{equation}

$\bullet \hspace{10pt} A=\diag \{a_{j}^{2}\}_{j=1}^{n}$, \hspace{5pt}$B=\diag \{b_{j}^{2}\}_{j=1}^{n}$ \hspace{10pt}($a_{j}>0$, $b_{j}>0$)\\
By Theorem 3.1, the heat kernel for the generalised Hermite operator has the form
\begin{equation}
\begin{split}
K(x,x^{0};t)&=(4\pi t)^{-\frac{n}{2}}\prod_{j=1}^{n}\left( \frac{2t b_{j}}{a_{j}\sinh(2t a_{j}b_{j})}\right) ^{\frac{1}{2}}\\
&\quad\times e^{-\frac{1}{4t}\left( \sum_{j=1}^{n}\frac{2t b_{j}}{a_{j}}\frac{\cosh(2t a_{j}b_{j})}{\sinh(2t a_{j}b_{j})}(x_{j}^{(1)})^{2}+\sum_{j=1}^{n}\frac{2t b_{j}}{a_{j}}\frac{\cosh(2t a_{j}b_{j})}{\sinh(2t a_{j}b_{j})}(x_{j}^{(0)})^{2}\right) }\\
&\quad\times e^{-\frac{1}{4t}\left(-2\sum_{j=1}^{n}\frac{2t b_{j}}{a_{j}}\frac{1}{\sinh(2t a_{j}b_{j})}x_{j}^{(0)}x_{j}^{(1)}\right) }.
\end{split}
\end{equation}

$\bullet \hspace{10pt} A>0$, \hspace{5pt}$B=\diag \{b_{j}^{2}\}_{j=1}^{n}$ \hspace{5pt}( $b_{j}>0$),\hspace{10pt} $A$ and $B$ are commutative.\\
With $D=4AB$ in Theorem 3.1, the heat kernel of the generalised Hermite operator is given by
\begin{equation}
\begin{split}
K(x,x^{0};t)&=(4\pi t)^{-\frac{n}{2}}\left(\frac{\det\psi\left(2t\sqrt{AB}\right)}{\det A} \right) ^{\frac{1}{2}}\times\\
&\quad e^{-\frac{1}{4t}\left(\left\langle \varphi\left(2t\sqrt{AB}\right)x, x \right\rangle _{A^{-1}}+\left\langle \varphi\left(2t\sqrt{AB}\right)x^{0}, x^{0} \right\rangle _{A^{-1}}-2\left\langle \psi\left(2t\sqrt{AB}\right)x, x^{0} \right\rangle _{A^{-1}}\right)}.
\end{split}
\end{equation}

\medskip
\begin{eg}
Ornstein-Uhlenbeck operator on weighted space
\end{eg}
Define Ornstein-Uhlenbeck operator 
$$
H_{OU}=-\div(A\nabla)+Bx\cdot\nabla
$$ 
with $A$ symmetric positive definite and $B$ any real matrix commutative with $A$. Ornstein-Uhlenbeck  
$$
H_{\phi}=-\div(A\nabla)+A\nabla\phi\cdot\nabla
$$ 
on Hilbert space $L^{2}(\R^{n}, e^{-\phi}dx)$ is unitarily equivalent to the Schr\"{o}dinger operator 
$$
H=-\div(A\nabla)+\frac{1}{4}\lvert\nabla\phi\rvert_{A}^{2}-\frac{1}{2}\div(A\nabla\phi)
$$
defined on Hilbert space $L^{2}(\R^{n}, dx)$: 
$$
H_{\phi}=THT^{-1}
$$ 
where $T$ is a multiplication operator defined by 
$$
Tu:=e^{\frac{\phi}{2}}u.
$$
Thus, 
$$
e^{-tH_{\phi}}=Te^{-tH}T^{-1}.
$$

Let $\phi$ take the form $\langle \tilde{B}x, x\rangle $ satisfying $A\nabla\phi=Bx$, then 
$$
\tilde{B}x+\tilde{B}^{t}x=\nabla\phi=A^{-1}Bx.
$$
Hence, 
$$
\phi=\frac{1}{2}\left\langle \left(\tilde{B}x+\tilde{B}^{t} \right)x,x \right\rangle=\frac{1}{2}\left\langle Bx, x \right\rangle _{A^{-1}},
$$
$$
\Hess(\phi)=A^{-1}B,
$$
$$
\div\left( A\nabla\phi\right)=\tr\left[ A\Hess(\phi)\right]=\tr(B).
$$
Consequently, 
$$
H=-\div\left( A\nabla\phi\right)+\frac{1}{4}\left\langle B^{t}A^{-1}Bx,x \right\rangle-\frac{1}{2}\tr(B).
$$
By Theorem 4.2 with $D=B^{t}B$, $f=g=0$ we have the heat kernel of $H$:
\begin{align*}
&\quad K(x,x^{0};t)\\
&=(4\pi t)^{-\frac{n}{2}}\left(\frac{\det\psi\left(t\sqrt{B^{t}B}\right)}{\det A} \right) ^{\frac{1}{2}}e^{\frac{t}{2}\tr(B)}\\
&\quad\times e^{-\frac{1}{4t}\left(\left\langle \varphi\left(t\sqrt{B^{t}B}\right)x, x \right\rangle _{A^{-1}}+\left\langle \varphi\left(t\sqrt{B^{t}B}\right)x^{0}, x^{0} \right\rangle _{A^{-1}}-2\left\langle \psi\left(t\sqrt{B^{t}B}\right)x, x^{0} \right\rangle _{A^{-1}}\right)}.
\end{align*}

For any $g\in L^{2}(\R^{n}, e^{-\phi}dx)$,
\begin{align*}
e^{-tH_{\phi}}g &= Te^{-tH}T^{-1}g\\
&= \int_{\R^{n}}e^{\frac{\phi(x)}{2}}K(x,y;t)e^{-\frac{\phi(y)}{2}}g(y)dy\\
&= \int_{\R^{n}}e^{\frac{\phi(x)}{2}}K(x,y;t)e^{\frac{\phi(y)}{2}}g(y)e^{-\phi}dy.
\end{align*}
Finally, the heat kernel of Ornstein-Uhlenbeck operator on weighted space 
$$
H_{OU}:\hspace{5pt} L^{2}(\R^{n}, e^{-\frac{\left\langle Bx, x \right\rangle _{A^{-1}}}{2}}dx) \hspace{5pt}\rightarrow \hspace{5pt} L^{2}(\R^{n}, e^{-\frac{\left\langle Bx, x \right\rangle _{A^{-1}}}{2}}dx)
$$
is given by 
\begin{equation}
\begin{split}
&\quad K_{OU}(x,x^{0};t)\\
&=e^{\frac{\phi(x)}{2}}K(x,x^{0};t)e^{\frac{\phi(x^{0})}{2}}\\
&=(4\pi t)^{-\frac{n}{2}}\left(\frac{\det\psi\left(t\sqrt{B^{t}B}\right)}{\det A}e^{t \tr(B)}\right) ^{\frac{1}{2}}\times\\
&\quad e^{-\frac{1}{4t}\left\lbrace \left\langle \left[ \varphi\left(t\sqrt{B^{t}B}\right)-tB\right]  x, x \right\rangle _{A^{-1}}+\left\langle\left[  \varphi\left(t\sqrt{B^{t}B}\right)-tB \right] x^{0}, x^{0} \right\rangle _{A^{-1}}-2\left\langle \psi\left(t\sqrt{B^{t}B}\right)x, x^{0} \right\rangle _{A^{-1}}\right\rbrace }.
\end{split}
\end{equation}

\bigskip
\noindent
Sheng-Ya Feng\\
Department of Mathematics\\
East China University of Science and Technology\\
Shanghai 200237, P.R. China\\
\medskip
E-mail address: s.y.feng@ecust.edu.cn
\end{document}